\input amstex.tex

\input amsppt.sty

\TagsAsMath

\magnification=1200

\hsize=5.0in\vsize=7.0in

\hoffset=0.2in\voffset=0cm

\nonstopmode

\document

\def\R{ \Bbb R}

\def\la{\langle}
\def\ra{\rangle}

\def\pd{\partial}

\input amstex.tex
\input amsppt.sty
\TagsAsMath \NoRunningHeads \magnification=1200
\hsize=5.0in\vsize=7.0in \hoffset=0.2in\voffset=0cm \nonstopmode

\document

\topmatter

\title{On asymptotic stability in energy space of  ground states of NLS in 2D}
\endtitle

\author
Scipio Cuccagna \,and \, Mirko Tarulli
\endauthor

\address
DISMI University of Modena and Reggio Emilia, via Amendola 2,
Padiglione Morselli, Reggio Emilia 42100 Italy\endaddress \email
cuccagna.scipio\@unimore.it \endemail

\address
Dipartimento di Matematica ``L. Tonelli'', University of Pisa, Largo Pontecorvo 5,
56127 Pisa, Italy\endaddress \email
tarulli\@mail.dm.unipi.it \endemail

\abstract We transpose  work by K.Yajima and by T.Mizumachi to prove
dispersive and smoothing
   estimates    for dispersive
solutions of the linearization at a ground state of a Nonlinear
Schr\"odinger equation (NLS) in 2D. As an application we extend to
dimension 2D a  result on asymptotic stability of ground states of
NLS proved in the literature for all dimensions different from 2.

\endabstract

\endtopmatter

 \head \S 1 Introduction \endhead

We consider even solutions  of a NLS
$$iu_t +\Delta u +\beta ( |u|^2 )u=0\, ,\, (t,x)
\in \Bbb R \times \Bbb R^2 \, , \,  u(0,x)=u_0(x) .\tag 1.1$$ We
assume:  {\item {(H1)}} $\beta (0)= 0$, $\beta\in C^\infty(\R,\R)$;
{\item {(H2)}}  there exists a $p_0\in(1,\infty)$ such that for
every $k=0,1$,
$$\left| \frac{d^k}{dv^k}\beta(v^2)\right|\lesssim
|v|^{p_0-k-1} \quad\text{if $|v|\ge 1$};$$ {\item {(H3)}} there
exists an open interval $\Cal{O}$ such that $
   \Delta u-\omega u+\beta(u^2)u=0
$ admits a $C^1$-family of ground states $\phi _ {\omega }(x)$ for
$\omega\in\Cal{O}$;

{\item {(H4)}} $ \frac d {d\omega } \| \phi _ {\omega
}\|^2_{L^2(\Bbb R )}>0$ for $\omega\in\Cal{O}$.

{\item {(H5)}} Let $L_+=-\Delta +\omega -\beta (\phi _\omega ^2
)-2\beta '(\phi _\omega ^2) \phi_\omega^2$ be the operator whose
domain is $H^2_{rad}(\Bbb R^2)$. We assume that $L_+$ has exactly
one negative eigenvalue.

\noindent By \cite{ShS} the $ \omega \to \phi _\omega \in H^1(\Bbb R
)$   is $C^2$ and by \cite{W1,GSS1-2} (H4-5) yields orbital
stability of the ground state $e^{i\omega t} \phi _ {\omega } (x)$.
Here we investigate asymptotic stability. We need some additional
hypotheses.

{\item {(H6)}} For any   $x\in \Bbb  {R} $, $u_0(x)=u_0(-x)$. That
is, the initial data $u_0$ of (1.1) are even.

Consider the Pauli matrices $\sigma _j$ and the linearization
$H_\omega$ given   by:
$$\aligned &\sigma _1=\left[ \matrix  0 &
1  \\
1 & 0
 \endmatrix \right] \, ,
\sigma _2=\left[ \matrix  0 &
i  \\
-i & 0
 \endmatrix \right] \, ,
\sigma _3=\left[ \matrix  1 &
0  \\
0 & -1
 \endmatrix \right]
\, ; \\&H_{\omega }=\sigma _3 \left [ - \Delta + \omega  -  \beta
(\phi ^2 _{\omega  }) - \beta ^\prime (\phi ^2 _{\omega  })\phi ^2
_{\omega   } \right ] +i  \beta ^\prime (\phi ^2 _{\omega  })\phi ^2
_{\omega  }
 \sigma _2    .
\endaligned \tag 1.2$$
Then we assume: {\item {(H7)}} Let $H_\omega$ be the linearized
operator around $e^{it\omega}\phi_\omega$, see (1.2).  $H_\omega$
has a positive simple eigenvalue $\lambda(\omega)$ for
$\omega\in\Cal{O}$. There exists an $N\in\Bbb N$ such that
$N\lambda(\omega)<\omega<(N+1)\lambda(\omega)$. {\item {(H8)}} The
Fermi Golden Rule (FGR) holds (see Hypothesis 4.2
  in Section 4).

{\item {(H9)}} The point spectrum of $H_\omega$ consists of $0$ and
$\pm\lambda(\omega)$.   The points $\pm \omega$ are not resonances.

Then we prove: \proclaim{Theorem 1.1}   Let $\omega_0\in \Cal {O}$
and $\phi_{\omega_0}(x)$ be a ground state  in a family of ground
states $\phi_{\omega } $. Let $u(t,x)$ be a solution to (1.1).
Assume (H1)--(H9). In particular assume the (FGR) in Hypothesis 4.2.
Then, if (1.1) is generic, there exist an $\epsilon_0>0$ and a $C>0$
such that for any $\epsilon \in (0,\epsilon _0)$ and for any $u_0$
with $ \|u_0-e^{i\gamma _0}\phi_ {\omega _0} \|_{H^1}<\epsilon  ,$
there exist $\omega_+\in \Cal {O}$, $\theta\in C^1(\Bbb R;\Bbb R)$
and $h_\infty \in H^1$ with $\| h_\infty\|
_{H^1}+|\omega_+-\omega_0|\le C \epsilon $ such that

$$
\lim_{t\to
+\infty}\|u(t,\cdot)-e^{i\theta(t)}\phi_{\omega_{+}}-e^{it\Delta
}h_\infty\|_{H^1}=0 .
$$
\endproclaim
Theorem  1.1  is the two dimensional version of  Theorem 1.1
\cite{CM}. The one dimensional version is in \cite{Cu3}. We recall
that  results of the sort discussed here were pioneered by Soffer \&
Weinstein \cite{SW1}, see also \cite{PW}, followed by Buslaev \&
Perelman \cite{BP1-2}, about 15 years ago. In this decade these
early works were followed by a number of results
\cite{BS,Cu1-2,GNT,M1,CZ,M2,P,RSS,SW1-3,TY1-3,Wd1}.  It was
heuristically understood that the rate of the leaking of energy from
the so called "internal modes" into radiation, is small and
decreasing when $N$ increases, producing technical difficulties in
the closure of the nonlinear estimates. For this reason prior to
Gang Zhou \& Sigal \cite{GS}, the literature treated only the case
when $N=1$ in (H6). \cite{GS} sheds light for $N>1$. The results in
\cite{GS} deal with all  spatial dimensions  different from 2 under
the so called Fermi Golden Rule (FGR) hypothesis. \cite{CM,Cu3}
strengthen \cite{GS} by considering initial data   in $H^1 $, by
showing that the (FGR) hypothesis is a consequence of what looks
generic condition, Hypothesis 4.2 below, if (H8) is assumed.
\cite{CM} treats also the case when there are many eigenvalues and
Hypothesis 4.2 is replaced by a more stringent hypothesis  which is
a natural generalization of the (FGR) hypothesis in \cite{GS}. The
same result  with many eigenvalues case can be proved also here and
in \cite{Cu3}, but we skip for simplicity the proof. We recall that
Mizumachi \cite{M1}, resp. \cite{M2}, extends to dimension 1, resp
2, the results in \cite{GNT} valid for  small solitons obtained by
 by bifurcation from ground states of a linear equation, while \cite{CZ} extends in 2D the
result in \cite{SW2}. \cite{Cu3} transposes
  \cite{M1}  to the case of large solitons, with the generalizations contained in \cite{CM}. Here
we consider the case of dimension 2.
  Thanks to the work by \cite{M2}, it is quite
clear how to transpose to dimension 2 the higher dimensional
arguments in \cite{CM}. The nonlinear arguments in \cite{CM} are not
sensitive to the dimension   except for the lack in 2D of the
endpoint Stricharz estimate. Mizumachi  \cite{M2} shows how to
replace it with   an appropriate smoothing estimate of Kato type.
The estimate and its proof are suggested by \cite{M2}.
 In order to complete the proof of Theorem
1.1 we need some dispersive estimates on the linearization $H_\omega
$ which in spatial dimension 2 are not yet proved in the literature.
The main technical task of this paper is the transposition to
$H_\omega $  of the proof of of $L^p$ boundedness of wave operators
of Schr\"odinger operators in dimension 2 due to Yajima \cite{Y2}.
We use the following notation. We set $
   H_0(\omega )=\sigma _3 (- \Delta + \omega  )
  $; given normed spaces $X$ and  $Y$ we denote by $B(X,Y)$ the
  space of operators from $X$ to $Y$ and given $L\in B(X,Y)$ we
  denote by $\| L \| _{X,Y}$  or by $\| L \| _{B(X,Y)}$ its norm.
We prove:

\proclaim{Proposition 1.2} Assume the hypotheses of Theorem 1.1.
The following limits are well defined isomorphism, inverse of each
other:
 $$\align  & W u= \lim _{t\to +\infty } e^{ itH_\omega } e^{- itH_0(\omega )
 }u  \text{ for any $u\in L^2$}  \\&
Z u= \lim _{t\to +\infty }  e^{ itH_0(\omega ) } e^{ -it H_\omega
 }  \text{ for any $u\in L^2_c(H_\omega )$ (defined in \S 2).}   \endalign $$
 For any $p\in (1,\infty )$  and any $k$ the restrictions of
$W$ and $Z$ to $L^2\cap W^{k,p}$ extend into operators such that for
$C(\omega )<\infty$ semicontinuous in $\omega$ $$\| W\| _{
W^{k,p}(\Bbb R^2), W^{k,p}  _c(H_\omega )}+\| Z\| _{
W^{k,p}_c(H_\omega ), W^{k,p}(\Bbb R^2) }<C(\omega ) $$ with
$W^{k,p}  _c(H_\omega ) $ the closure in $W^{k,p} (\Bbb R^2)$ of $
W^{k,p}(\Bbb R^2)\cap L^2  _c(H_\omega ) $.
\endproclaim

We will set  $L^{2,s}$ and $H^{m,s}$

$$\|u\|_{L^{2,s}}=\|\langle  x\rangle ^su\|_{L^2(\R^2)}
\quad\text{and}\quad\|u\|_{H^{m,s}}=\| \langle x\rangle
^su\|_{H^m(\R^2)},$$
 where $m\in \Bbb N$, $s\in \R$ and $\langle x\rangle =(1+|x|^2)^{1/2}$.
For $f(x)$ and $g(x)$ column vectors, their inner product is $
\langle f ,
  g\rangle =\int _{\Bbb R^2}{^tf(x)}\cdot  \overline{{ {g(x)}}} dx $. The
  adjoint $H^\ast$ is defined by $
\langle Hf ,
  g\rangle =
\langle f ,
  H^\ast g\rangle .$ Given an operator $H$, its resolvent is
  $R_H(z)=(H-z)^{-1}.$ We will write $R_0(z)=(-\Delta -z)^{-1}.$
  We write
  $  \| g(t,x) \| _{L^p_tL^q_x}= \|  \| g(t,x) \| _{ L^q_x} \| _{L^p_t
  } $ and $  \| g(t,x) \| _{L^p_tL^{2,s}_x}= \|  \| g(t,x) \| _{ L^{2,s}_x} \| _{L^p_t
  } $

\head \S 2 Linearization, modulation and set up \endhead

We will use  the  following classical result,  \cite{We1,GSS1-2},
see also \cite{Cu3}:

\proclaim{Theorem 2.1}  Suppose that $e^{i\omega t} \phi _ {\omega }
(x)$ satisfies (H4).   Then $\exists \, \epsilon >0$ and a
$A_0(\omega )>0$ such that for any  $\| u(0,x) - \phi _ {\omega } \|
_{H^1}<\epsilon $ we have for the corresponding solution $ \inf  \{
\|  u(t,x) -e^{i \gamma }\phi _ {\omega } (x-x_0) \| _{H^1(x\in \Bbb
R ^2)} : \gamma \in \Bbb R \, \& \, x_0 \, \in \Bbb R  \}  <
A_0(\omega ) \epsilon . $ \endproclaim      We can write the ansatz
$  u(t,x) = e^{i \Theta (t)} (\phi _{\omega (t)} (x)+ r(t,x)) \, ,
\, \Theta (t)= \int _0^t\omega (s) ds +\gamma (t). $ Inserting the
ansatz into the equation we get
$$\aligned &
  i r_t  =
 -  r _{xx}+\omega (t) r-
\beta ( \phi _{\omega (t)} ^2 )r -\beta ^\prime ( \phi _{\omega (t)}
^2 )\phi _{\omega (t)} ^2 r \\&-
 \beta ^\prime ( \phi _{\omega (t)} ^2 )
\phi _{\omega (t)} ^2  \overline{  r }+ \dot \gamma (t) \phi
_{\omega (t)} -i\dot \omega (t)
\partial _\omega \phi   _{\omega (t)}
+ \dot \gamma (t) r
 +
  O(r^2).\endaligned
$$
 We set
$^tR= (r,\bar r) $,  $^t\Phi = ( \phi _{\omega } , \phi _{\omega } )
$
 and we rewrite the above equation as
$$i  R _t =H _{\omega}   R +\sigma _3 \dot \gamma   R
+\sigma _3 \dot \gamma \Phi - i \dot \omega \partial _\omega \Phi
+O(R^2).\tag 2.1
$$
 Set
  $H_0(\omega )=\sigma _3(-d^2/dx^2 +\omega )$
and $V(\omega )=H_\omega -  H_0(\omega ).$ The essential spectrum is
$$\sigma _e =\sigma _e (H_\omega )
=\sigma _e (H_0(\omega ) ) =(-\infty , -\omega ] \cup [ \omega ,
+\infty ) .$$  0 is an isolated eigenvalue. Given an operator $L$ we
set $N_g(L)= \cup _{j\ge 1} \ker(L^j)$. \cite{We2}
 implies that,
 if $\{ \cdot \}$ means span,
$N_g(H^\ast _\omega )=\{ \Phi , \sigma _3\partial _\omega  \Phi   \}
$.  $\lambda  (\omega )$ has corresponding
   real eigenvector  $\xi  (\omega )$, which can be normalized  so that $\langle \xi   , \sigma _3 \xi   \rangle
 =1$.    $\sigma _1\xi (\omega ) $ generates $\ker (H  _\omega
+\lambda (\omega ))$ .  The function  $(\omega , x )\in \Cal O
\times \Bbb R  \to \xi  (\omega , x)$ is $C^2$; $|\xi  (\omega , x)|
< c e^{-a|x|}$ for fixed $c>0$ and $a>0$ if $\omega \in K \subset
\Cal O$, $K$ compact. $\xi   (\omega , x)$ is even in $x$ since by
assumption we are restricting ourselves in the category of such
functions. We have  the
 $H _{\omega}$   invariant
Jordan block decomposition
$$\align & L^2= N_g(H _{\omega} ) \oplus \big (
\oplus _{j,\pm  }\ker (  H _{\omega}\mp \lambda  (\omega )) \big )
\oplus L_c^2(H _{\omega})=N_g(H _{\omega})\oplus
 N_g^\perp (H _{\omega}^\ast  )
\endalign
$$
where we set $L_c^2(H _{\omega} )=\left \{   N_g(H _{\omega} ^\ast )
\oplus \oplus _{ \pm  }\ker (H _{\omega} ^\ast \mp \lambda (\omega
))\right \} ^{\perp}.$ We can impose
$$    R (t) =(z  \xi + \bar z  \sigma _1 \xi ) + f(t)   \in
\big [ \sum _{ \pm  } \ker (H _{\omega (t)}\mp \lambda  (\omega
(t)))\big ] \oplus L_c^2(H _{\omega (t)})  .\tag 2.2  $$ The
following claim admits an elementary proof which we skip:

\proclaim{Lemma 2.2} There is a Taylor expansion at $R=0$ of the
nonlinearity $O(R^2)$ in (2.1) with $R_{m,n}(\omega  ,x) $ and
$A_{m,n}(\omega ,x ) $ real vectors  and matrices rapidly decreasing
in $x$: $ O(R^2)=$
$$\aligned &  \sum _{ 2\le  m+n \le 2N+1} R_{m,n}(\omega ) z^m  \bar z^n+
\sum _{1\le  m + n \le  N} z^m  \bar z^n A_{m,n}(\omega ) f+
O(f^2+|z|^{2N+2}) .\endaligned  $$
\endproclaim
In terms of the frame in (2.2) and the expansion in Lemma 2.2, (2.1)
becomes
$$\aligned &if_t=\left ( H _{\omega (t)}+\sigma _3 \dot \gamma \right )f + \sigma _3
\dot \gamma \Phi (\omega )- i \dot \omega \partial _\omega \Phi (t)
+
  (z \lambda  (\omega ) -i\dot z ) \xi  (\omega ) \\&
- (\bar z \lambda  (\omega )+i\dot {\bar z }) \sigma _1\xi (\omega )
  +\sigma _3 \dot \gamma (z  \xi + \bar z  \sigma _1 \xi )
-i \dot \omega (z  \partial _\omega \xi + \bar z  \sigma _1
\partial _\omega \xi ) \\& + \sum _{ 2\le  m+n \le 2N+1} z^m \bar
z^nR_{m,n}(\omega ) + \sum _{1\le  m + n \le N} z^m \bar z^n
A_{m,n}(\omega ) f+\\& + O(f^2)+ O_{loc}(|z  ^{2N +2}|)
\endaligned \tag 2.3
$$
where by $O_{loc}$ we mean that the there is a factor $\chi (x)$
rapidly decaying   to 0   as $|x|\to \infty $. By taking inner
product of the equation with generators of $N_g(H_\omega ^\ast )$
and $\ker (H_\omega ^\ast -\lambda  )$ we obtain modulation and
discrete modes equations:

$$\aligned & i\dot \omega \frac{d\| \phi _\omega \| _2^2}{d\omega  }
=\langle \sigma _3 \dot \gamma (z \xi + \bar z  \sigma _1 \xi ) -i
\dot \omega (z \partial _\omega \xi + \bar z  \sigma _1
\partial _\omega \xi )  +   \sum _{    m+n =2}^{ 2N+1}
 z^m \bar z^nR_{m,n}(\omega )\\& +  \big ( \sigma _3 \dot \gamma +i\dot \omega
\partial _\omega P_c+ \sum _{  m + n =1}^{ N} z^m \bar z^n
A_{m,n}(\omega )  \big ) f   + O(f^2)+ O_{loc}(|z ^{2N +2}|), \Phi
\rangle \\& \dot \gamma \frac{d\| \phi _\omega \| _2^2}{d\omega  }
=\langle \text{ same as above }, \sigma _3
\partial _\omega \Phi \rangle \\& i\dot z -\lambda (\omega )z =\langle
\text{ same as above }, \sigma _3 \xi  \rangle .
\endaligned \tag 2.4$$

\head \S 3 Spacetime estimates  for $H_\omega $
\endhead

We need analogues of Lemmas 2.1-3 and Corollary 2.1 in \cite{M2}. We
call admissible  all pairs $(p,q)$ with $ 1/p=1/2-1/q$ and $2\le q
<\infty$.  We set $(p',q')=(p/(p-1),q/(q-1))$. In the lemmas below
we assume that the $H_\omega$ of the form (1.2) for which hypotheses
(H3-5), (H7) and (H9) hold.

\proclaim{Lemma 3.1 (Strichartz estimate)} There exists a positive
number $C=C(\omega )$ upper semicontinuous in $\omega $ such that
for any $k\in [0,2]$:
  {\item {(a)}}
 for any $f\in
L^2_c( {\omega })$ and any admissible  all pairs $(p,q)$,
$$\|e^{-itH_{\omega }} f\|_{L_t^pW_x^{k,q } }\le C\|f\|_{H^k}.$$
 {\item {(b)}}
  for any $g(t,x)\in
S(\R^2)$ and any couple of admissible pairs $(p_{1},q_{1})$
$(p_{2},q_{2})$ we have
$$
\|\int_{0}^te^{-i(t-s)H_{\omega }} P_c( {\omega
})g(s,\cdot)ds\|_{L_t^{p_1}W_x^{k,q_1 } } \le
C\|g\|_{L_t^{p_2'}W_x^{k,q_2' } }.
$$\endproclaim
Lemma 3.1 follows immediately from Proposition 1.2 since $W$ and $Z$
intertwine $e^{-itH_{\omega }}P_c(H_\omega )$  and $e^{-itH_{0 }} $.

\medskip

\proclaim{Lemma 3.2} Let $s>1$. $\exists$   $C=C(\omega )$ upper
semicontinuous in $\omega $ such that:
   {\item {(a)}}
  for any $f\in S(\R ^2 )$,
$$\align & \|  e^{-itH_{\omega }}P_c( {\omega })f\|_{  L^2_tL_x^{2,-s}} \le
 C\|f\|_{L^2} ;
\endalign $$
  {\item {(b)}}
  for any $g(t,x)\in
 {S}(\R^3)$
$$ \left\|\int _\Bbb R e^{itH_{\omega }}P_c( {\omega })g(t,\cdot)dt\right\|_{L^2_x} \le
C\| g\|_{ L_t^2L_x^{2,s}}.
$$\endproclaim

Notice that (b) follows from (a) by duality.

\proclaim{Lemma 3.3} Let $s>1$.  $\exists$   $C=C(\omega )$ as above
such that $\forall$ $g(t,x)\in {S}(\R^3)$ and $t\in\R$:
$$\align &  \left\|  \int_0^t e^{-i(t-s)H_{\omega
}}P_c( {\omega })g(s,\cdot)ds\right\|_{  L^2_tL_x^{2,-s}} \le C\|
g\|_{  L^2_tL_x^{2, s}} .\endalign
$$
\endproclaim

As a corollary from Christ and Kiselev \cite{CK}, Lemmas 3.2 and
 3.3 imply:

\proclaim{Lemma 3.4} Let $(p,q)$ be an admissible pair and let
$s>1$. $\exists$   $C=C(\omega )$ as above such that $\forall$
$g(t,x)\in {S}(\R^3)$ and $t\in\R$:
$$
\left\|\int_0^t e^{-i(t-s)H_{\omega }}P( {\omega })g(s,\cdot)ds
\right\|_{L_t^pL_x^q} \le C\|g\|_{  L^2_tL_x^{2, s}}.$$
\endproclaim

\proclaim{Lemma 3.5} Consider the diagonal matrices $
E_+=\text{diag}(1 , 0)$ $E_-=\text{diag}(0 , 1).$ Set $P_\pm (\omega
)=Z(\omega ) E_\pm W(\omega )$ with $Z(\omega )$ and $W(\omega )$
the wave operators associated to $H_\omega $. Then  we have for
$u\in L^2_c(H_\omega )$
$$\aligned &P_+(\omega )u =\lim _{\epsilon \to 0^+}
 \frac 1{2\pi i}
\lim _{M \to +\infty} \int _\omega ^M \left [ R_{H_\omega }(\lambda
+i\epsilon )- R_{H_\omega }(\lambda -i\epsilon )
 \right ] ud\lambda \\&
P_-(\omega )u =\lim _{\epsilon \to 0^+}
 \frac 1{2\pi i}
\lim _{M \to +\infty} \int _{-M }^{-\omega } \left [ R_{H_\omega
}(\lambda +i\epsilon )- R_{H_\omega }(\lambda -i\epsilon ) \right ]
ud\lambda
\endaligned \tag 1
$$
and  for any $s_1$ and $s_2$ and for $C=C (s_1,s_2,\omega  )$ upper
semicontinuous in $\omega $, we have $$  \| (P_+(\omega
 )-P_-(\omega  )-P_c(\omega  )\sigma _3) f\|  _{L^{2,s_1} }\le
C  \|       f\|  _{L^{2,s_2} }. \tag 2$$
\endproclaim
{\it Proof.} Formulas (1) hold  with $P_{\pm}(\omega )$ replaced by
$E_\pm $ and $H_\omega $ replaced by $H_0$ and for any $u\in
L^2(\Bbb R^2)$. Applying $W(\omega )$ we get (1) for $H_\omega $.
Estimate (2) follows by the proof of inequality (3)   in Lemma 5.12
\cite{Cu3} which is valid for all dimensions.

\head \S 4 Proof of Theorem 1.1 \endhead

   We restate Theorem 1.1 in a more precise form:

\proclaim{Theorem 4.1} Under the assumptions of Theorem 1.1   we can
express
$$u(t,x)=e^{i \Theta (t)}   (\phi _{\omega (t)} (x)+  \sum _{j=1}^{2N} p_j(z,\bar z)  A _j (x,\omega (t))+
 h(t,x)  )$$ with $p_j(z,\bar z)=O(z )$ near 0,
with  $\lim _{t\to +\infty}\omega (t)$ convergent, with $ |A _j
(x,\omega (t)) |\le C   e^{-a|x|}$ for fixed $C>0$ and $a>0$, $\lim
_{t\to +\infty} z(t  )  =0,$ and   for fixed $C>0$
$$\| z(t)\|_{ L_t ^{2N+2}}^{N+1} +\| h(t,x)\| _{L^\infty  _tH ^{1 }_x \cap  L^3 _tW ^{1,6}_x }<C \epsilon .\tag 1$$
   Furthermore, there exists
$h_\infty\in H^1(\Bbb R ,\Bbb C )$ such that
$$\lim_{t\to\infty}\|  e^{i\int _0^t\omega (s) ds +i\gamma (t)}h(t)-
e^{it  \Delta   }h_\infty\|_{H^1}=0. \tag 2$$
\endproclaim

The proof of Theorem 4.1  consists in a normal forms expansion and
in the closure of some nonlinear estimates. The normal forms
expansion is exactly the same of \cite{CM,Cu3}, in turn adaptations
of \cite{GS}.

\head \S 4.1 Normal form expansion\endhead

We  repeat \cite{CM}.
  We   pick $k=1,2,...N$ and
set $f=f_k$ for $k=1$. The other  $f_k$ are defined below. In the
ODE's  there will be error terms    of the form
$$E_{ODE}(k)=   O( |z |^{2 N+2 }
 )+O(  z^{N+1}   f _{k}  )  +O(f^2_{k})+
O(\beta (|f _{k}|^2)f _{k} ).$$ In the PDE's there will be error
terms   of the form
$$E_{PDE}(k)=   O_{loc}( |z |^{N+2}
 )  +O_{loc}(  z f _{k}  )+O_{loc}(f^2_{k})+
O(\beta (|f _{k}|^2)f _{k} ).$$ In the right hand sides of the
equations  (2.3-4) we substitute $\dot \gamma $ and $\dot \omega $
using the modulation equations. We repeat the procedure a sufficient
number of times until we can write for $k=1$ and $ f_1=f$

$$\aligned   i\dot \omega \frac{d\| \phi _\omega \| _2^2}{d\omega  }
=&\langle   \sum _{   m +n =2}^{ 2N+1 } z^m \bar z^n\Lambda
_{m,n}^{(k)}(\omega ) +   \sum _{   m +n =1}^{ N } z^m \bar z^n
A_{m,n}^{ (k) }(\omega ) f_{k}  +E_{ODE}(k) , \Phi (\omega ) \rangle
  \\  i\dot z -\lambda  z =&\langle
\text{ same as above }, \sigma _3 \xi (\omega )  \rangle
\\ i\partial _t f_{k}=&\left ( H _{\omega}  +\sigma _3 \dot \gamma
\right  )f_{k} + E_{PDE}(k)+ \sum _{k+1\le  m +n \le N +1}z^m \bar
z^n R_{m,n}^{(k)}(\omega ),
\endaligned
$$
with  $A_{m,n}^{(k)} $, $R_{m,n}^{(k)}     $  and $\Lambda
_{m,n}^{(k)} (\omega , x) $ real
 exponentially  decreasing  to 0 for $|x|\to \infty$ and continuous in $(\omega , x) $.
Exploiting $|(m-n)\lambda (\omega )|<\omega $ for $  m+n \le N$,
$m\ge 0$, $n\ge 0$,
 we define inductively $f_k$ with $k\le N$ by

 $$f_{k-1}=-\sum _{  m+n=k}z^m \bar z^n R_{H_{\omega }}((m-n)\lambda (\omega ) ) R_{m,n}^{(k-1)}   (\omega
 )
+f_k.$$ Notice that if $ R_{m,n}^{(k-1)}   (\omega
 ,x)$ is real
 exponentially  decreasing  to 0 for $|x|\to \infty$, the same is
 true for $R_{H_{\omega }}((m-n)\lambda (\omega ) ) R_{m,n}^{(k-1)}   (\omega
 )$ by $|(m-n)\lambda (\omega )|<\omega $. By induction $f_k$ solves the above equation with
the above notifications. Now we manipulate the   equation for $f_N$.
We fix $\omega _1=\omega (0)$. We write

$$\aligned &i\partial _t P_c(\omega _1)f_{N}-\left \{ H_{\omega _1}   + (\dot \gamma +\omega -\omega _1)
(P_+(\omega _1)-P_-(\omega _1))\right \} P_c(\omega _1)f_{N} =\\&
+P_c(\omega _1)\widetilde{E}_{PDE}(N) +\sum _{  m+n= N+1} z^m \bar
z^n P_c(\omega _1)R_{m,n}^{(N)} (\omega _1)
\endaligned \tag 4.1$$ where we split $P_c(\omega _1)=P_+(\omega _1)+P_-(\omega
_1)$ with $P_\pm (\omega _1)$, see Lemma 3.5, where $P_+(\omega _1)$
are the projections in $\sigma _c(H_{\omega _1}) \cap \{ \lambda :
\pm \lambda \ge \omega _1 \} $ and with
$$\aligned & \widetilde{E}_{PDE}(N)= E_{PDE}(N)  +  \sum _{  m+n= N+1} z^m \bar
z^n  \left ( R_{m,n}^{(N)} (\omega  )-R_{m,n}^{(N)} (\omega _1)
\right ) +  \varphi (t,x  ) f_N\\& \varphi (t,x  ):=
 (\dot \gamma +\omega -\omega _1) \left (P_c(\omega _1)\sigma _3-  (P_+(\omega _1)-P_-(\omega _1))  \right ) f_N +
   \left (  V  (\omega ) -
 V  (\omega _1) \right )   f_N
\\& +(\dot \gamma +\omega -\omega _1) \left ( P_c(\omega  )-  P_c(\omega _1)\right )\sigma _3f_N.
  \endaligned \tag 4.2
$$
By Lemma 3.5   for $C_N(\omega _1)$ upper semicontinuous in $\omega
_0$,   $\forall $ $N$ we have
$$ \|  \langle x \rangle ^{N}  (P_+(\omega
_1)-P_-(\omega _1)-P_c(\omega _1)\sigma _3) f\|  _{L^2_x}\le
C_N(\omega _1) \|  \langle x \rangle ^{-N}    f\|  _{L^2_x}.$$  The
term  $\varphi (t,x  )$ in (4.2) can be treated as a small cutoff
function. We write
$$\aligned &f_{N}=-\sum _{  m+n=N+1} z^m \bar z^nR_{H_{\omega _1}}((m-n)  \lambda (\omega _1) +i0 )
 P_c(\omega
_1)R_{m,n}^{(N)}   (\omega _1)  +  f_{N+1}.\endaligned \tag 4.4$$
Then
$$\aligned &i\partial _t P_c(\omega _1) f_{N+1}=\left ( H_{\omega _1}    + (\dot \gamma +\omega -\omega _1)
(P_+(\omega _1) -P_-(\omega _1) )\right ) P_c(\omega _1) f_{N+1} +
\\& +  \sum _{\pm } O(\epsilon |z|^{N+1}
)R_{H_{\omega _1}} (\pm (N+1) \lambda (\omega _1 )+i0 ) R_\pm
(\omega _1) +P_c(\omega _1)\widehat{ {E}}_{PDE}(N)   \endaligned
\tag 4.5$$ with $R_+ =R_{ N+1,0}^{(N)} $ and $R_- =R_{ 0, N+1
}^{(N)} $ and $\widehat{ {E}}_{PDE}(N)=\widetilde{ {E}}_{PDE}(N)+
O_{loc}(\epsilon z^{N+1})$, where we have  used that $( \omega
-\omega _0)=O(\epsilon )$ by Theorem 2.1. Notice that $R_{H_{\omega
_0}} (\pm (N+1) \lambda (\omega _0)+i0 ) R_\pm (\omega _0)\in
L^\infty$ do not decay spatially.
  In the ODE's with $k=N$, by the standard theory of normal forms
and following the idea in Proposition 4.1 \cite{BS}, see \cite{CM}
for details, it is possible to introduce new unknowns
$$\aligned &\widetilde{\omega }=\omega +q(\omega , z,\bar z) +\sum _{1\le  m
+n \le N }z^m  \bar z^n\langle f_N,\alpha _{mn}(\omega )\rangle ,\\&
\widetilde{z } =z  +p (\omega , z,\bar z) +\sum _{1\le  m +n \le N
}z^m \bar z^n\langle f_N,\beta _{mn}(\omega )\rangle ,\endaligned
\tag 4.6$$ with $p (\omega , z,\bar z)=\sum p_{ m,n}(\omega )
z^m\bar z^n$ and $q(z,\bar z)=\sum q_{ m,n}(\omega ) z^m\bar z^n$
polynomials in $(z,\bar z)$ with real coefficients and $O(|z|^2)$
near 0, such that we get
$$\aligned & i  \dot {\widetilde{\omega}} =  \langle  {E}_{PDE}(N) , \Phi \rangle   \\ &
 i\dot {\widetilde{z} }-\lambda  (\omega )\widetilde{z}  =
 \sum _{ 1\le  m \le N} a_{ m  }(\omega )|\widetilde{z}^{m}| ^{2 }
 \widetilde{z} +\langle  E_{ODE}(N) , \sigma _3 \xi   \rangle +\\& +
  \overline{\widetilde{z}}^N \langle  A_{0,N}^{(N)}(\omega ) f _{N} , \sigma _3
\xi   \rangle  .
\endaligned \tag 4.7$$
with $   a_{ m  }(\omega )  $ real. Next step is to substitute $
f_{N} $ using (4.4). After eliminating by a new change of variables
$\widetilde{z}= \widehat{{z}}+p(\omega
,\widehat{z},\overline{\widehat{z}})$ the resonant terms, with
$p(\omega ,\widehat{z},\overline{\widehat{z}})=\sum \widehat{p}_{
m,n}(\omega ) z^m\bar z^n$ a polynomial in $(z,\bar z)$ with real
coefficients $O(|z|^2)$ near 0, we get

$$\aligned & i  \dot {\widehat{\omega}} =  \langle  {E}_{PDE}(N) , \Phi \rangle   \\ &
 i\dot {\widehat{z} }-\lambda  (\omega )\widehat{z}  =
 \sum _{ 1\le m\le N} \widehat{a}_{ m  }(\omega )|\widetilde{z}^{m}| ^{2 }
 \widehat{z} +\langle  E_{ODE}(N) , \sigma _3 \xi   \rangle   -\\&
 -
  |\widehat{z} ^N|^2\widehat{z}  \langle  \widehat{A}_{0,N}^{(N)}(\omega )
  R_{H_{\omega _0}}((N+1) \lambda  (\omega _1) +i0)P_c(\omega
_0)R_{N+1 ,0}^{(N)}(\omega _1)     ,\sigma _3\xi   \rangle
\\& +  \overline{\widehat{z}} ^N \langle  \widehat{A}_{0,N}^{(N)}(\omega ) f
_{N+1} , \sigma _3 \xi   \rangle
\endaligned \tag 4.8$$
with $\widehat{a}_{ m  }$, $\widehat{A}_{0,N}^{(N)}$ and
$R_{N+1,0}^{(N)}$ real. By $\frac{1}{x-i0}=PV\frac{1}{x}+i\pi \delta
_0(x)$ and by an elementary use of the wave operators, we can denote
by $\Gamma (\omega , \omega _0)$ the quantity
$$\aligned &\Gamma  (\omega ,
\omega _1)= \Im \left (\langle  \widehat{A}_{0,N}^{(N)}(\omega )
  R_{H_{\omega  _1}}((N+1) \lambda (\omega _1) +i0)P_c(\omega
_1)R_{N+1,0}^{(N)}(\omega _1) \sigma _3 \xi   (\omega )\rangle
\right )\\& =\pi
  \langle  \widehat{A}_{0,N}^{(N)}(\omega )
  \delta (H_{\omega  _1}- (N+1) \lambda (\omega _1)  )P_c(\omega
_1)R_{N+1,0}^{(N)}(\omega _1) \sigma _3 \xi   (\omega )\rangle .
 \endaligned $$
Now we assume the following: \proclaim{Hypothesis 4.2} There is a
fixed constant $\Gamma >0$ such that $|\Gamma  (\omega , \omega
)|>\Gamma .$
\endproclaim
By continuity and by Hypothesis 4.2 we can assume $|\Gamma (\omega ,
\omega _1 )|>\Gamma /2.$ Then we write

$$\aligned &  \frac{d}{dt} \frac{|{\widehat{z} }|^2}{2} =    -\Gamma  (\omega , \omega  _1)
  |z|  ^{2N+2} +\Im \left (
\langle  \widehat{A}_{0,N}^{(N)}(\omega )f _{N+1} , \sigma _3 \xi
(\omega )\rangle \overline{\widehat{z}} ^{N+1}   \right ) \\&+ \Im
\left ( \langle E_{ODE}(N) , \sigma _3 \xi  (\omega ) \rangle
\overline{\widehat{z}} \right ) .
\endaligned \tag 4.9
$$

\head \S 4.2 Nonlinear estimates \endhead

 By an elementary
continuation argument, the following a priori estimates imply
inequality (1) in Theorem 4.1, so to prove (1) we focus on:

\proclaim{Lemma 4.3} There are fixed  constants $C_0$ and $C_1$ and
$\epsilon _0>0$  such that   for any $0<\epsilon \le \epsilon _0$
  if we have
$$\| \widehat{ z} \| _{L^{2N+2 }_t}^{N+1} \le 2C_0\epsilon \quad \&
 \quad \| f_N \| _{L^\infty  _tH ^{1 }_x \cap L^3 _tW ^{1,6}_x \cap L^{\frac{
2{p_0}}{{p_0}-1}}_tW^{1,2{p_0}}_x \cap L^2_tH^{1,-s}} \le
2C_1\epsilon \tag 4.10$$ then we obtain the improved inequalities
$$  \align&
\| f_N \| _{L^\infty  _tH ^{1 }_x \cap L^3 _tW ^{1,6}_x \cap
L^{\frac{ 2{p_0}}{{p_0}-1}}_tW^{1,2{p_0}}_x \cap L^2_tH^{1,-s}} \le
C_1\epsilon  , \tag 4.11
\\& \| \widehat{z} \| _{L^{2N+2 }_t}^{N+1} \le C_0\epsilon .\tag
4.12 \endalign$$
\endproclaim
{\it Proof}. Set $\ell (t):=\gamma +\omega -\omega _1$.  First of
all, we have:

\proclaim{Lemma 4.4} Let $g (0,x)\in H^1_x\cap L^2_c(\omega _1)$ and
let $\omega (t)$ be  a continuous function. Consider     $  i g _t=
\left \{ H _{\omega _1}   + \ell (t) (P_+(\omega _0)-P_-(\omega
_0))\right \}g+P_c(\omega _1) F.$  Then for a fixed $C=C(\omega
_1,s)$ upper semicontinuous in $\omega _1$ and $s>1$ we have
$$\aligned & \| g \| _{L^\infty _tH ^{1 }_x \cap L^3 _tW ^{1,6}_x
\cap L^{\frac{ 2{p_0}}{{p_0}-1}}_tW^{1,2{p_0}}_x }
 \le  C   ( \| g(0,x)\| _{H^1}+\| F\|  _{   L^1
_tH^1 _x  +   L_t^2H_x^{1,s} }).\endaligned$$
\endproclaim
Lemma 4.4 follows easily from Lemmas 3.1-4 and $ P_\pm (\omega
_1)g(t)=$
$$=e^{-it H_{\omega _1}}e^{- i\int _0^t \ell  (\tau ) d\tau }P_\pm
(\omega _1)g(0)  -i\int _0^{t }e^{-i(t-s) H_{\omega _1}} e^{\pm
i\int _s^t \ell (\tau ) d\tau }P_{ \pm}(\omega _1)
  F  (s) ds
$$

 \proclaim{Lemma 4.5} Consider equation (4.1) for
 $f_N$ and assume (4.10). Then we can split $\widetilde{E}_{PDE}(N)=
 X+ O(f_N^3)+O( f _{N}^{p_0} )$ such that
 $ \|  X\| _{    L_t^2H_x^{1,M} }  \lesssim \epsilon ^2
 $ for any fixed $M$ and  $ \|O(f_N^3)+O( f _{N}^{p_0} )\| _{ L^1
_tH^1 _x}   \lesssim \epsilon ^3.  $
\endproclaim
{\it Proof of Lemma 4.5.} In the error terms for $k=N$ at the
beginning of \S 4.1 we can write $\widetilde{E}_{PDE}(N) = $
$$O(\epsilon ) \psi (x) f_N +O_{loc}( |z |^{N+2}
 )  +    O_{loc}(  z f _{N}  )+O_{loc}(f^2_{N})+
O(   f _{N}^3 )+O( f _{N}^{p_0} ) $$ with $\psi (x)$ a rapidly
decreasing function, $p_0$ the exponent in (H2) and with $O( f
_{N}^{p_0} )$ relevant only for $p_0>3$.  Denoting $X$ the sum of
all terms except the last one,   setting $f=f_N$, by (4.10)  we
have:
 :\medskip {\item{(1)}} $ \| O(\epsilon ) \psi (x) f \| _{
L_t^2H_x^{1,M} } \lesssim \epsilon \|   f  \| _{
L_t^2H_x^{1,-M}}\lesssim \epsilon ^2;$
\medskip {\item{(2)}}$ \| O_{loc}(  z f    ) \| _{ L_t^2H_x^{1,M} } \lesssim \| z\| _{\infty} \|  f  \| _{
L_t^2H_x^{1,-M}}\lesssim \epsilon ^2 ;$
\medskip
{\item{(3)}} $ \| O_{loc}(    f  ^2  ) \| _{ L_t^2H_x^{1,M} }
\lesssim   \|   f  \| _{ L_t^2H_x^{1,-M}}^2\lesssim \epsilon ^2.$

This yields $ \| \la x \ra ^{  M} X\| _{   H_x^1 L_t^2 }  \lesssim
\epsilon ^2
 $. To bound the remaining term  observe:

{\item{(4)}} $ \|     |f  |^2  f  \| _{  L_t^1H^1_x } \lesssim \left
\|  \| f \| _{W^{1,6}_x}    \| f \| ^2_{L^{ 6}_x} \right \| _{L^1_t
} \le \| f \| _{L^3_tW^{1,6}_x}^3 \lesssim \epsilon ^3;$

{\item{(5)}} $ \|    O(  f  ^{{p_0}} ) \| _{  L_t^1H^1_x } \lesssim
\left \|  \| f \| _{W^{1,2{p_0}}_x}    \| f \| ^{{p_0}-1}_{L^{
2{p_0}}_x} \right \| _{L^1_t } \le \| f \| _{L^{\frac{
2{p_0}}{{p_0}-1}}_tW^{1,2{p_0}}_x}
 \| f \| ^{{p_0}-1}_{L^{  2{p_0}  \frac{{p_0}-1}{{p_0}+1}  }_tW^{1,2{p_0}}_x} \lesssim \epsilon ^{p_0},$
where in the last step we use $ \| f \| _{L^{  2{p_0}
\frac{{p_0}-1}{{p_0}+1} }_tW^{1,2{p_0}}_x}\lesssim \| f \| ^{\alpha
}_{L^{\frac{ 2{p_0}}{{p_0}-1}}_tL^{ 2{p_0}}_x} \|
   f      \| _{  L_{t }^\infty H^1_x } ^{1-\alpha }$ for some
   $0<\alpha <1$ by ${p_0}>3$, interpolation and Sobolev embedding.

\bigskip

 {\it Proof of (4.11).} Recall that $f_N$ satisfies equation (4.1)
whose right hand side is $P_c(\omega _1)\widetilde{E}_{PDE}(N) +
O_{loc}(    z  ^{N+1}  )$. In addition to Lemma 4.5 we have the
estimate $ \| O_{loc}(    z  ^{N+1}  ) \| _{ L_t^2H_x^{1,M} }
\lesssim \| z \| _{ L_t^{2N+1} }^{N+1} \lesssim 2C_0\epsilon  .$  So
by Lemmas 3.1-4, for some fixed $c_2$ we get schematically
$$\| f_N \| _{L^\infty  _tH ^{1 }_x \cap L^3 _tW ^{1,6}_x \cap L^{\frac{
2{p_0}}{{p_0}-1}}_tW^{1,2{p_0}}_x }  \le 2c_2C_0\epsilon +\epsilon
+O(\epsilon ^2)$$ where $\epsilon$ comes from initial data,
$O(\epsilon ^2)$ from all the nonlinear terms save for
 the $ R_{m,n}^{(N)}   (\omega _0) z^m  \bar z^n $ terms which contribute the  $2c_2C_0\epsilon
 $. Let now $f_N= g+h$ with
$$\aligned & i g _t=
\left \{ H _{\omega _1}   + \ell (t) (P_+(\omega _1)-P_-(\omega
_1))\right \}g+X +O_{loc}(    z  ^{N+1}  )\, , \quad  g(0)=f_N(0)\\&
i h _t= \left \{ H _{\omega _1}   + \ell (t) (P_+(\omega
_1)-P_-(\omega _1))\right \}h+O(f_N^3)+O( f _{N}^{p_0} ) \, , \quad
h(0)=0\endaligned
$$ in the notation of Lemma 4.5. Then, by Lemmas 3.2 and 3.3 and by
the estimates in Lemma 4.5 we get $\|   g \| _{ L^2_tH_x^{1,-s}}
\lesssim 2C_0\epsilon +O(\epsilon ^2)+c_0\epsilon $ for a fixed
$c_0$. Finally,
$$\aligned &\int _0^\infty  \|e^{-i(t-s) H_{\omega _1}}  e^{\pm i \int _s^t\ell (\tau ) d\tau }( O(f_N^3)+O( f _{N}^{p_0}
))(s)\| _{L^2_tH^{1,-s}} \\& \lesssim \int _0^\infty  \| (
O(f_N^3)+O( f _{N}^{p_0} ))(s)\| _{ H^{1 }}\lesssim \epsilon ^3.
\endaligned $$
 So if we set  $C_1\approx 2C_0+c_0+1$ we obtain (4.11). We
need to bound $C_0$.

\medskip

{\it Proof of (4.12).}
 We first need: \proclaim{Lemma 4.6} We can
decompose $f _{N+1}= h_1+h_2+h_3 +h_4$ with   for  a fixed large
$M>0$: {\item {(1)}}       $\|   h _1\| _{L^2_tL _x^{2,M}} \le
O(\epsilon ^{2})   ;$ {\item {(2)}}   $\|
  h _2\| _{L^2_tL _x^{2,M}}\le O(\epsilon ^{2}) ;$

{\item {(3)}}   $\|   h _3\| _{L^2_tL _x^{2,M}} \le O(\epsilon
^{2});$ {\item {(4)}}   $\|   h _4\| _{L^2_tL _x^{2,M}} \le c(\omega
_1)\epsilon  $ for a fixed $c(\omega _1) $ upper semicontinuous in
$\omega _1$.
\endproclaim
{\it Proof of Lemma 4.6}. We set $$\aligned & i\partial _t h_1=\left
( H_{\omega _1}   + \ell (t) (P_+-P_-)\right ) h_1\\& h_1(0)=\sum _{
m+n=N+1} R_{H_{\omega _1}}((m-n)\lambda (\omega _1)+i0)
R_{m,n}^{(N)} (\omega _1)z^m(0)  \bar z^n (0) .\endaligned $$ We get
$\|   h _1\| _{L^2_tL _x^{2,-M}} \le c(\omega _1) |z(0)|^2 \sum \|
  R_{m,n}^{(N)} (\omega _1)\|
_{ L _x^{2,M} }=O(\epsilon ^2) $ by  the following lemma:

\proclaim{Lemma 4.7} There is a fixed $s_0$ such that for $s>s_0$,
$$\aligned &
\|   e^{-iH _\omega t}R_{H_\omega} (\Lambda +i0) P_c(\omega )
\varphi \| _{L^2_tL^{2,-s}_x} < C_s(\Lambda , \omega ) \|
  \varphi  (x) \| _{ L^{2, s}_x}   \\&  \left \|  \int _0^t e^{-iH _\omega (t-\tau )}R_{H_\omega} (\Lambda +i0) P_c(\omega )
g(\tau ) d\tau \right \| _{L^2_tL^{2,-s}_x} < C_s(\Lambda , \omega )
\|
  g(t,x) \| _{L^2_t L^{2, s}_x}  \endaligned \tag 4.13
$$
with   $ C_s(\Lambda , \omega )$ upper semicontinuous in $\omega $
and in $\Lambda >\omega $.
\endproclaim
Let us assume Lemma 4.7  for the moment, for the proof see \S 9. We
set $h_2(0)=0$ and

$$\aligned &i\partial _t h_{2}=\left ( H_{\omega _1}   + \ell (t)  (P_+-P_-)\right ) h_{2}
+\\& + O(\epsilon z^{N+1}) R_{H_{\omega _1}}((N+1)\lambda (\omega
_1)+i0 ) R_{N+1,0}^{(N)} (\omega   _0)
\\& +
O(\epsilon z^{N+1}) R_{H_{\omega _1}}(-(N+1)\lambda (\omega _1)+i0 )
R_{0, N+1}^{(N)} (\omega   _1) .\endaligned $$
 Then we have
$h_{2}=h_{21 }+h_{22}$ with $h_{2j}=\sum _\pm h_{2j\pm }$ with $h_{2
1\pm } (t)=$

$$\int _0^te^{-iH_{\omega _1}(t-s)}e^{\pm i\int _s^t \ell(\tau )
d\tau }P_{ \pm}
 O(\epsilon z^{N+1}) R_{H_{\omega _1}}((N+1)\lambda (\omega _1)+i0 )  R_{N+1,0}^{(N)}   (\omega   _1) ds $$
and  $h_{22\pm }$ defined similarly but with $R_{H_{\omega
_0}}(-(N+1)\lambda (\omega _1)+i0 ) R_{0,N+1 }^{(N)}$ . Now by
(4.13) we get
$$\|h_{2j\pm } (t)   \| _{L^2_tL^{2, -M}_x} \le C \epsilon \| z \|  _{L^{ 2N+2}
_{t }} ^{N+1} $$ and so $  \|h_{2 } (t)   \| _{L^2_tL^{2, -M}_x}
   = O(\epsilon ^2).$
Let   $h_3(0) =0$   and

$$\aligned &i\partial _t P_c(\omega _1) h_3=\left ( H_{\omega _1}   + \ell (t)
(P_+(\omega _1) -P_-(\omega _1) )\right ) P_c(\omega _1) h_3 +
 P_c(\omega _1) \widetilde{E}_{PDE}(N). \endaligned
 $$
 Then by the argument in the proof of (4.11) we get claim (3).
Finally let   $h_4(0) =f_N(0)$   and

$$\aligned &i\partial _t P_c(\omega _1) h_4=\left ( H_{\omega _1}   + \ell (t)
(P_+(\omega _1) -P_-(\omega _1) )\right ) P_c(\omega _1) h_4  .
\endaligned
 $$
Then by  Lemma 3.2   $ \| \langle x\rangle ^{-M} h _4\| _{L^2_{tx}}
\lesssim   \| f_N(0)\| _{L^2_x}\le c(\omega _1)\epsilon $ we get
(4).
\bigskip

{\it Continuation of proof of Lemma 4.3}. We integrate (4.9) in
time. Then by Theorem 2.1  and by Lemma 4.4 we get, for $A_0$   an
upper bound  of the    constants $A_0(\omega )$ of Theorem 2.1,

$$ \| \widehat{z}\| _{L^{2N+2}_t}^{2N+2}\le  A_0\epsilon ^2+\epsilon
\| \widehat{z}\| _{L^{2N+2}_t}^{N+1} + o(\epsilon ^2).$$ Then we can
pick    $C_0=( A_0+1)$ and this proves that (4.10) implies (4.12).
Furthermore $\widehat{z}(t)\to 0$ by
$\frac{d}{dt}\widehat{z}(t)=O(\epsilon ).$

\bigskip

As in \cite{CM,Cu3} in the above argument we did not use the sign of
$\Gamma (\omega , \omega _0)$. With the same argument in
\cite{CM,Cu3}
 one can prove

 \proclaim{Corollary 4.8} If Hypothesis 4.2 holds, then $\Gamma (\omega , \omega
 )>\Gamma $.
 \endproclaim

\bigskip
The proof that, for $^tf_N(t)=(h(t),\overline{h}(t))$, $h(t)$ is
asymptotically free for $t\to \infty$, is similar to the analogous
one in \cite{CM} and we skip it.

\head \S 5 Limiting absorption principle and $L^2$ theory for
$H_\omega $
\endhead
In sections \S 5- \S 7 we prove Proposition 1.2.  We start
emphasizing two consequences of hypothesis (H9), in particular (b)
clarifies the absence of resonance at $\pm \omega$:

{\item {(a)}} $H_\omega$ has no eigenvalues in $[\omega , +\infty
)\cup (-\infty ,-\omega ]$;

{\item {(b)}} if   $g\in W^{2,\infty}(\Bbb R^2, \Bbb C^2)$ satisfies
$H_\omega g=\omega g$   or $H_\omega g=-\omega g$ then $g=0$.

Because of the fact that  $H_\omega$ is not a symmetric operator, we
need some preparatory work to show that in fact $H_\omega$ is
diagonalizable in the continuous spectrum. This work is done in \S 5
which ends with a formula for the wave operator $W$ which is the
basis to develop in \S 6-7 a transposition of the work of Yajima
\cite{Y2}.

We first need a preliminary on Schr\"odinger operators. We will
denote by  $q(x)$  a real valued
 function with:  $q(x)\ge
0$ with $q(x)>0$ at some points;   $q(x)\in C^\infty _0(\Bbb R^2)$.
We set $h_q=-\Delta +q(x)$. Then we have:

\proclaim{Lemma 5.1} Let $\Bbb C_+ =\{ z\in \Bbb C : \, \Im z>0\}$.
Suppose $q(x)=0$ for $r\ge r_0>0$. Then we have the following facts.

{\item {(1)}} There exists $ s_0>0$ and $C_0>0$ such that  for $s\ge
s_0$,  $R_{h_q} (z)$ extends into a function   $ z\to R_{h_q}^{+}(z)
$ which is in $ (L^\infty \cap C^0)( \overline{\Bbb C_+}, B
(L^{2,s}, L^{2,-s} ))$.

{\item {(2)}} For any $n_0\in \Bbb N$  there exists $ s_0>0$ such
that   for any $a_0>0$ there is a choice of $C>0$ such that for
$n\le n_0$
$$\left  \| \frac{d^n}{dz^n} R_{h_q}^+(z) : L^{2,s}(\Bbb R^2)\to L^{2,-s}(\Bbb
R^2)\right \| \le C_0 \langle z \rangle ^{-\frac{1}{2} (1+n)} \text{
$\forall$ $z\in \overline{\Bbb C_+} \cap \{ z:  |z|  \ge a_0\} .$}
  $$

{\item {(3)}} The same argument can be repeated for $\Bbb C_- =\{
z\in \Bbb C : \, \Im z<0\}$ and $R_{h_q}^-(z)$.
\endproclaim
Claim (2) follows from \cite{Ag} and \cite{JK} and claim (3) follows
along the lines of the previous two claims. In view of (2), it is
enough to prove (1) for $z\approx 0.$
 For $\zeta
=re^{i\theta}$ with $\theta \in (-\pi , \pi )$ let $\sqrt{\zeta}
=\sqrt{r}e^{i\theta /2}$. With this convention  for $z\not \in
[0,\infty )$ for $R_0(z)=(-\Delta -z)^{-1}$ we have  $$
R_0(z)=\frac{1}{2\pi}K_0(\sqrt{-z}|x|)\ast=\frac{i}{4}H^{+}_0(
i\sqrt{-z}|x|)\ast =-\frac{i}{4}H^{-}_0( -i\sqrt{-z}|x|)\ast $$ for
the Macdonald function $K_0$ and the Hankel functions $H^{\pm }_0$.
We set $G_0=- \frac{1}{2\pi} \log |x|\ast$,  $P_0 f =\int _{\Bbb
R^2}fdx$. We have  for $M(z)=(1+\sqrt{q}R_{ 0} (z)\sqrt{q})
 $ the identity
$$R_{h_q} (z)= R_{ 0} (z)-R_{ 0} (z)\sqrt{q} M^{-1}
(z)\sqrt{q}R_{ 0} (z) .\tag 4$$     From the expansion at 0 in
$\Bbb C_+$ of $H^{+}_0$  and by the argument in Lemma 5 \cite{Sc} we
have in $B(L^{2,s} ,L^{2,-s} ),$ for $s $ sufficiently large,
$$R_{0}  (z) =c (z) P_0 - G_0 +O(-z\log \sqrt{-z})  \quad
 c (z)= \frac{i}{4}- \frac{\gamma }{2\pi} -\frac{1}{ 2\pi}\log(\sqrt{-z}/2).
 \tag 5$$
Consider the  projections in $L^2(\Bbb R^2)$, $P=\sqrt{q}\langle
\cdot , \sqrt{q}\rangle /\| q\| _{L^1}$ and $Q=1-P$. Let
$T=1+\sqrt{q}G_{ 0}  \sqrt{q}$. Then $QTQ$ is invertible in $QL^2(\Bbb
R^2)$. Denote its inverse in $QL^2(\Bbb R^2)$ by $D_0=(QTQ)^{-1}$.
Consider the operator in $L^2=PL^2\oplus QL^2$ defined by
$$S= \left [ \matrix   P  &  & -PTQD_0Q \\ -QD_0 QTP & & QD_0QTPTQD_0Q \endmatrix    \right
]
$$
and $h(z)=\| q\| _{L^1}c(z)+\text{trace}(PTP-PTQD_0QTP)$. Then by
\cite{Sc}
$$\aligned & R_{h_q} (z)= R_{ 0} (z)-h^{-1}(z)R_{ 0} (z)\sqrt{q}S \sqrt{q}R_{ 0}
(z) \\& -R_{ 0} (z)\sqrt{q}QD_0Q \sqrt{q}R_{ 0} (z)-R_{ 0}
(z)\sqrt{q}O(-z\log \sqrt{-z})\sqrt{q}R_{ 0} (z) .\endaligned \tag 6$$ By
direct computation
$$\aligned & h^{-1}(z)R_{ 0} (z)\sqrt{q}S \sqrt{q}R_{ 0}
(z) =  \frac{c^2(z)}{h(z)}\langle \cdot , 1 \rangle \sqrt{q}S
\sqrt{q}\langle \cdot , 1 \rangle + \frac{c (z)}{h(z)}\langle \cdot
, 1 \rangle \sqrt{q}S \sqrt{q}G_0 +\\& +  \frac{c (z)}{h(z)}G_0
\sqrt{q}S \sqrt{q}\langle \cdot , 1 \rangle +  \frac{c (z)}{h(z)}G_0
\sqrt{q}S \sqrt{q} G_0 +O(-z\log \sqrt{-z}),
\endaligned  $$
where all terms, except the first on the right hand side, 
admit continuous extension in $\overline{\Bbb C} _+$ at 0. We have $
\langle \cdot , 1 \rangle \sqrt{q}S \sqrt{q}\langle \cdot , 1
\rangle = \| q\| _{L^1} P_0  $ and so by (5)

$$ R_{ 0} (z)-\frac{c^2(z)}{h(z)} \| q\| _{L^1}
P_0$$ admits continuous extension in $\overline{\Bbb C} _+$ at 0. By
direct computation
$$ \aligned & R_{ 0} (z)\sqrt{q}QD_0Q \sqrt{q}R_{ 0} (z)=G_0\sqrt{q}QD_0Q
\sqrt{q}G_0+O(-z\log \sqrt{-z})
\endaligned
$$
  admits continuous extension in $\overline{\Bbb C} _+$ at 0. So $R_{ h_q} (z)
$ admits continuous extension in $\overline{\Bbb C} _+$ at 0, and so
on all $\overline{\Bbb C} _+$.

\bigskip
A   consequence of Lemma 5.1 is the $h_q$ smoothness in the sense of
Kato \cite{Ka} of multiplication operators involving rapidly
decreasing functions $\psi$:

 \proclaim{Lemma 5.2} Let $\psi
(x)\in L^\infty (\Bbb R^2)\cap L^{2,s}(\Bbb R^2)$ for $s\gg 1$   and
$q$ as in Lemma 5.1. Then the multiplication operator $\psi $ is
$h_q$ smooth, that is, for a fixed $C>0$
$$ \int _{\Bbb R}\|  \psi R_{h_q} (\lambda +i\varepsilon )
u\|_2^2d\lambda < C\|   u\|_2^2 \text{ for all $u\in L^2(\Bbb R^2)$
and $\varepsilon \neq 0$} .$$
\endproclaim
This follows from one of the characterizations of  $H$   smoothness
in the   case $H$ is selfadjoint, see Theorem 5.1 \cite{Ka},
specifically from the fact that by Lemma 5.1 we have that for $\psi
_1, \psi _1\in L^\infty  \cap L^{2,s} $ for $s\gg 1$ there is  a
number $C>0$ such that for all $z\not \in \Bbb R$ we have $  \| \psi
_1R_{h_q} (z ) \psi _2 \| _{L^2, L^2}<C.$

\bigskip

We consider now  $
   H_q=\sigma _3 (- \Delta + q+\omega  )
  $ and consider our linearization  $H_\omega $.
Write  $H_\omega =H_q+( V_\omega -\sigma_3q)$, and factorize
$V_\omega -\sigma_3q =B^\ast A$ with $A,B$ smooth $|\partial
_x^\beta A(x)|+|\partial _x^\beta B(x)| < C e^{-\alpha |x|}$
$\forall \, x,$ for some $\alpha ,C >0$ and for $|\beta |\le N_0$,
$N_0$ sufficiently large.
     We have $\sigma
_1H_q=-H_q\sigma _1$, $\sigma _1H_\omega =-H_\omega \sigma _1$. We
choose the factorization   $B^\ast A$ so that $\sigma _1B^\ast
=-B^\ast \sigma _1$, $\sigma _1A = A \sigma _1$. By these equalities
 $\sigma _1R _{ H_q} (z)=-R_{ H_q}(-z)\sigma _1$ and $\sigma _1R_{H  _\omega }(z)=-R_{H  _\omega }(-z)\sigma
 _1$, so in some of the estimates below it is enough to consider
 $z\in \overline{\Bbb C_{++}}$ with $ \Bbb C_{++}=\{ z:\Im z>0, \, \Re z>0
 \} .$

\proclaim{Lemma 5.3} For $ z \in \overline{\Bbb C_{+ }}$ the
function $ R_{H_q}^+(z) $ is well defined and satisfies the
following properties:

{\item{ (1)}} There exists $ s_0>0$ and $C_0>0$ such that  for $s\ge
s_0$  the function   $ z\to R_{H_q}^+(z) $   is in $( L^\infty \cap
C^0)( \overline{\Bbb C_+}, B (L^{2,s}, L^{2,-s} ))$.

{\item {(2)}} For any $n_0\in \Bbb N$  there exists $ s_0>0$ such
that   for any $a_0>0$ there is a choice of $C>0$ such that for
$n\le n_0$ and  $\forall$ $z\in \overline{\Bbb C_+} \cap \{ z:
\text{dist} (z,\pm \omega ) \ge a_0\} ,$
$$\left  \| \frac{d^n}{dz^n} R_{H_q}^+(z) : L^{2,s}(\Bbb R^2)\to L^{2,-s}(\Bbb
R^2)\right \| \le C_0 \langle z \rangle ^{-\frac{1}{2} (1+n)} .
  $$

{\item {(3)}} For any $\psi (x)\in L^\infty (\Bbb R^2)\cap
L^{2,s}(\Bbb R^2)$ for $s\gg 1$  the multiplication operator $\psi $
is $H_q$ smooth, that is, for a fixed $C>0$
$$ \int _{\Bbb R}\|  \psi R_{H_q} (\lambda +i\varepsilon )
u\|_2^2d\lambda < C\|   u\|_2^2 \text{ for all $u\in L^2(\Bbb R^2)$
and $\varepsilon \neq 0$} .$$

{\item {(4)}} Analogous statements hold for $ z \in \overline{\Bbb
C_{- }}$ and the function $ R_{H_q}^-(z) $.

\endproclaim

Lemma 5.3 is a trivial consequence of Lemmas 5.1-2. The properties
in Lemma 5.4 are partially inherited by $H_\omega$. Let $
{Q}_q^+(z)= {A}R_{H_q}^+(z) {B}^\ast $. Then for $z\in \Bbb C_+$

\proclaim{Lemma 5.4} Fix an exponentially decreasing bounded
function $\psi  $. For  $ z \in  {\Bbb C_{+ }}$ the function $
AR_{H_\omega } (z) \psi $ extends into a function $ AR_{H_\omega
}^{+} (z) \psi $ for  $ z \in  \overline{ \Bbb C_{+ }} \backslash
\sigma _d(H_\omega )$ with the following properties:

{\item{ (1)}} $\forall$ $a_0>0$ $\exists$   $C_0>0$ such that for $
X_{a_0}=\overline{\Bbb C_+}\cap \{ z: \text{dist} (z,\sigma
_d(H_\omega ) ) \ge a_0\} $

$$
AR_{H_\omega }^{+} (z) \psi \in  ( L^\infty \cap C^0) ( X_{a_0} , B
(L^{2 }, L^{2 } )) $$

{\item {(2)}} For any $n_0\in \Bbb N$  there exists $ s_0>0$ such
that   for any $a_0>0$ there is a choice of $C>0$ such that for
$n\le n_0$ and  $\forall$ $z\in X_{a_0} \cap \{ z: \text{dist}
(z,\pm \omega ) \ge a_0\} ,$
$$\left  \| \frac{d^n}{dz^n} AR_{H_\omega }^{+} (z) \psi  : L^{2 }(\Bbb R^2)\to L^{2 }(\Bbb
R^2)\right \| \le C_0 \langle z \rangle ^{-\frac{1}{2} (1+n)} .
  $$

{\item {(3)}} There is  a constant $C>0$ such that
$$\int \|   {A} R _{H_\omega }(  \lambda +i \varepsilon )  u\| ^2_2 d\lambda \le
C \| u\| ^2_2  \text{ for all $ u \in L_c^2(H _\omega )   $ and
$\varepsilon \neq 0$} .$$

{\item {(4)}} Analogous statements hold for $ z \in \overline{\Bbb
C_{- }}$ and the function $ R_{H_\omega }^-(z) $.

\endproclaim
{\it Proof.} Let us write $ {Q}_q^{+}(z)=  {A}R_{H_q}^+(z) {B}^\ast
$ and for $z\in \Bbb C _{+}$
$$ {A}R_{H_\omega} (z  ) =
 ( 1+  {Q}_q^+(z )  )^{-1} A R_{H_q}(z  ) .\tag 5
$$
By Lemma 5.3 we have  $\lim _{z\to \infty } \|
 {Q}_q^+(z)\| _{ L^2 , L^2}   =0  $. By analytic
Fredholm theory $1+  {Q}_q^+(z )$ is not invertible  only at the
  $z\in \overline{\Bbb C _{+}}$ where $\ker ( 1+  {Q}_q^+(z )
)\neq 0$. This set has 0 measure in $\Bbb R$. By Lemma 2.4
\cite{CPV} if at some $z\neq \pm \omega $ we have $\ker ( 1+
{Q}_q^+(z) )\neq 0$, then $z$ is an eigenvalue. By hypothesis there
are no eigenvalues in $\sigma _e(H_\omega ).$ Hence we get claim
(2).

\proclaim{Lemma 5.5} If $\ker ( 1+ Q_q^+ (\omega ) )\neq 0$ then
there exists $g\in W^{2,\infty}(\Bbb R^2)$ with $g \neq 0$ such that
$H_\omega g=\omega g$
\endproclaim
Let us assume Lemma 5.5. By hypothesis such $g$ does not exist. This
yields (1). By  (5), claim (4) Lemma 5.4 and Neumann expansion we
get (4). Next,  apply (5) to $u \in L_c(H _\omega )$. $
{A}R_{H_\omega} (z  )u$ is an analytic function in $z$ with values
in $L^2(\Bbb R^2)$ for $z$ near any isolated eigenvalue $z_0$ of
$H_\omega $ because the natural projection of $u $ in $N_g(H_\omega
-z_0)$ is 0. Away from isolated eigenvalues of $H_\omega $,
 $ ( 1+  {Q}_q^+(z )  )^{-1}$ is uniformly bounded. Hence (3) in
 Lemma 5.3 implies (3) in Lemma 5.4.

 {\it Proof of Lemma 5.5.} Let $0\neq \widetilde{g}\in \ker ( 1+ Q_q^+ (\omega )
 )$. Then $$B^\ast \widetilde{g}+(V_\omega -q)R_{H_q} (\omega  )B^\ast
 \widetilde{g}=0.$$ Set $g=R_{H_q} (\omega  )B^\ast
 \widetilde{g}$. Then $Ag=- \widetilde{g}$ and so $g\neq 0$. By $g+R_{H_q} (\omega
 ) (V_\omega -q) g=0 $ we have $g\in H^2 _{ loc}(\Bbb
 R^2)$ and   $H_\omega g=\omega g$. We want now to show that
 $g\in L^\infty (\Bbb R^2)$, contrary to the hypotheses.
We have $^tg=(g_1,g_2)$ with $ g_2=(\Delta -q-2\omega ) ^{-1}
(B^\ast
 \widetilde{g})_2$, where $B^\ast
 \widetilde{g}\in L^{2,s}(\Bbb R^2)$ for any $s$, so $g_2\in H^2(\Bbb R^2)$. We have
 $ g_1= R_{h_q} ^{+}(0)(B^\ast
 \widetilde{g})_1$  with $g_1\in L ^{2,-s}(\Bbb R^2)$ for
 sufficiently large $s$. We split $L^{2,\pm s}= L^{2,\pm s}_r\oplus \left ( L^{2,\mp
s}_r\right ) ^{\perp}$ where $L^{2,\pm s}_r$ are the radial
functions and we are considering the standard pairing $L^{2,s
}\times L^{2,-s}\to \Bbb C$ given by $\int _{\Bbb R^2}f(x) g(x) dx$.
We decompose  $g_1=g_{1r}+g_{1nr}$ with $g_{1r}\in L ^{2,-s} _{r}$
 and $g_{1nr}\in ( L ^{2, s} _{r})^\perp .$ In  $( L ^{2, -s} _{r})^\perp  \to ( L ^{2, s} _{r})^\perp  $  we have
 $ R_{h_q}^+(0)= G_0-G_0q (1+QG_0qQ)
^{-1}G_0$   with $Q=1-P$, for $P=P_0q_0$, $q_0=c_0^{-1}q$, $c_0=\int
_{\Bbb R^2}q dx$, $P_0u=\int _{\Bbb R^2}u dx$.   Then
$$g_{1nr}=G_0(B^\ast
 \widetilde{g})_{ 1nr}-G_0q (1+QG_0qQ)
^{-1}G_0 (B^\ast
 \widetilde{g})_{ 1nr}$$
and by   asymptotic expansion for $|x|\to \infty$ we conclude that
for some constants
$$\partial _x^\alpha  \left (g_{1nr}-a-\frac{b_1x_1+b_2x_2}{|x|^2}\right
) =O(|x|^{-1-\alpha -\epsilon})  $$ for some $\epsilon >0$. Finally
we look ar $\widetilde{g} _{ 1r}.$ We can consider   solutions $
\phi (r ) $ and $ \psi (r ) $   of $h_qu=0$    with: $ \phi (0 ) =1$
and $ \phi _{r} (0 ) =0$;
   $\psi (r_0 ) =1$ and  $|\psi (r  )  |$  bounded for $r\ge
 r_0$, $\psi (r_0 )\approx c\log r $ with $c\neq 0$ for $r\to 0$.  In terms of these two functions the kernel of
 $R_{h_q}^+(0)$ in $L^2((0,\infty ),dr)$ is
$$\aligned  R_{h_q}^+(0)(r_1,r_2) =  \frac {\phi (r_1 ) \psi (r_2 )  }
{W(r_2)} \, \text{if} \, r_1<r_2 \text{ or } =\frac {\phi (r_2 )
\psi (r_1 )  } {W(r_2)}  \, \text{if} \, r_1>r_2,
\endaligned
$$
with $W(r)=[ \phi (\cdot ),  \psi (\cdot  )] (r)=c / r$ for some
$c\neq 0$. We have $g_{1 r}(r)=$
$$\aligned &  = c^{-1} \psi (r  ) \int _0^r\phi (s )(B^\ast
 \widetilde{g})_{ 1 r} (s) \, s\, ds  +c^{-1} \phi (r  ) \int _ r^{+\infty }\psi (s )(B^\ast
 \widetilde{g})_{ 1 r} (s)\, s\, ds.
\endaligned $$
Then for $r\ge r_0$, $| g_{1 r}(r)|\le$

$$\aligned &     |c^{-1} \psi (r  )| \int _0^r|\phi (t )(B^\ast
 \widetilde{g})_{ 1 r} (t)| \, t\, dt  +|c^{-1} \phi (r  )| \int _ r^{+\infty }|\psi (t )(B^\ast
 \widetilde{g})_{ 1 r} (t)|\, t\, dt \\& \lesssim \| \log \langle x \rangle  \|  _{L^{2,-s}(\Bbb R^2)}
 \| B^\ast
 \widetilde{g}    \|  _{L^{2,s}(\Bbb R^2)} + \log (2+r)  \| B^\ast
 \widetilde{g}    \|  _{L^{2,s}(   \{ x\in \Bbb R^2 :|x|\ge r \})} =O(1)  .
\endaligned $$
Then we conclude that we have a nonzero $g\in H^2 _{ loc}(\Bbb
 R^2)\cap  L^\infty (\Bbb R^2)$ such that $H_\omega g=\omega g$. But
 this is contrary to the nonresonance hypothesis.

 \bigskip Analogous to Lemma 5.4 is:

\proclaim{Lemma 5.6} Fix an exponentially decreasing bounded
function $\psi  $. For  $ z \in  {\Bbb C_{+ }}$ the function $ B
R_{H_\omega ^\ast  } (z) \psi $ extends into a function $B
R_{H_\omega ^\ast  }^{+} (z) \psi $ for  $ z \in  \overline{ \Bbb
C_{+ }} \backslash \sigma _d(H_\omega )$ with the following
properties:

{\item{ (1)}} For any $a_0>0$ there exists   $C_0>0$ such that $ B
R_{H_\omega ^\ast  }^{+} (z) \psi \in  L^\infty ( X_{a_0} , B (L^{2
}, L^{2 } )) $ where $ X_{a_0}=\overline{\Bbb C_+}\cap \{ z:
\text{dist} (z,\sigma _d(H_\omega ) ) \ge a_0\}.$

{\item {(2)}} For any $n_0\in \Bbb N$  there exists $ s_0>0$ such
that   for any $a_0>0$ there is a choice of $C>0$ such that for
$n\le n_0$ and  $\forall$ $z\in X_{a_0} \cap \{ z: \text{dist}
(z,\pm \omega ) \ge a_0\} ,$
$$\left  \| \frac{d^n}{dz^n} B
R_{H_\omega ^\ast  }^{+} (z) \psi  : L^{2 }(\Bbb R^2)\to L^{2 }(\Bbb
R^2)\right \| \le C_0 \langle z \rangle ^{-\frac{1}{2} (1+n)} .
  $$

{\item {(3)}} There is  a constant $C>0$ such that
$$\int \|   B
R_{H_\omega ^\ast  }(  \lambda +i \varepsilon )  u\| ^2_2 d\lambda
\le C \| u\| ^2_2  \text{ for all $ u \in L_c^2(H _\omega ^\ast)   $
and $\varepsilon \neq 0$} .$$

{\item {(4)}} Analogous statements hold for $ z \in \overline{\Bbb
C_{- }}$ and the function $ R_{H_\omega  ^\ast }^-(z) $.

\endproclaim

From \S 2 \cite{Ka}  we conclude: \proclaim{Lemma  5.7}   There  are
isomorphisms
 $\widetilde{W}\colon L^2\to L^2_c(H_\omega )  $ and $\widetilde{Z}\colon L^2_c(H_\omega )  \to  L^2  $,
inverses of each other, defined as follows: for $u\in  L^2$,  $v\in
L^2_c(H_\omega ^\ast )   $,
$$\aligned & \langle \widetilde{W} u,v\rangle =
\langle  u,v\rangle  +\lim _{\epsilon \to 0^+ } \frac 1{2\pi i }
\int _{-\infty}^{+\infty} \langle  {A}  R_{H_q}(  \lambda +i
\epsilon )  u, {B} R_{H^\ast _\omega }   (\lambda +i\epsilon )
v\rangle d\lambda ;
\endaligned$$
for $u\in   L^2_c(H_\omega ) $,  $v\in  L^2 $,
$$\aligned & \langle \widetilde{Z} u,v\rangle =
\langle  u,v\rangle  +\lim _{\epsilon \to 0^+ } \frac 1{2\pi i }
\int _{-\infty}^{+\infty} \langle  {A} R_{H  _\omega } (\lambda
+i\epsilon ) u, {B} R_{H_q}(\lambda +i\epsilon )v\rangle d\lambda .
\endaligned$$
We have $H_\omega \widetilde{W}=\widetilde{W}H_q$ and $
{H}_q\widetilde{Z}=\widetilde{Z}H_\omega $, $ e^{ itH_\omega }
\widetilde{W}=\widetilde{W} e^{ itH_q }$ and $ e^{ itH_q
}\widetilde{Z}=\widetilde{Z}e^{ itH_\omega }P_c(H_\omega ). $ The
operators $\widetilde{W}$ and $\widetilde{Z}$ depend continuously on
$\widetilde{A}$ and $\widetilde{B}^\ast $ and can be expressed as
$$\align  & \widetilde{W }u= \lim _{t\to +\infty } e^{ itH_\omega } e^{
-itH_q
 }u  \text{ for any $u\in L^2$}  \\&
\widetilde{Z} u= \lim _{t\to +\infty }  e^{ itH_q } e^{- it H_\omega
 }  \text{ for any $u\in L^2(H_\omega )$} . \endalign $$
\endproclaim
In particular we remark:\proclaim{Lemma 5.8} We have for $C(\omega
)$ upper semicontinuous in $\omega $ and
$$ \| e^{-itH_\omega } g \| _2\le C(\omega ) \|  g\| _2 \text{ for any $g\in
L^2_{c}(H_\omega )$}.$$

\endproclaim

Having proved that $e^{-itH_\omega }  P_c(H_\omega )$ are bounded in
$L^2$, we want to relate $H_\omega$ to $
   H_0=\sigma _3 (- \Delta + \omega  )
  $ .
Write  $H=H_0+V_\omega $, $V_\omega =B^\ast A$.  We have $\sigma
_1H_0=-H_0\sigma _1$, $\sigma _1H_\omega =-H_\omega \sigma _1$. We
choose the factorization of $V_\omega$ so that $\sigma _1B^\ast
=B^\ast \sigma _1$, $\sigma _1A =-A \sigma _1$. By these equalities
 $\sigma _1R_{H_0}(z)=-R_{H_0}(-z)\sigma _1$ and $\sigma _1R_{H  _\omega }(z)=-R_{H  _\omega }(-z)\sigma
 _1$.
We have the following result about existence and completeness of
wave operators:
 \proclaim{ Lemma 5.9} The following limits are well defined:
 $$\align  & W u= \lim _{t\to +\infty } e^{ itH_\omega } e^{- itH_0
 }u  \text{ for any $u\in L^2$} \tag 1\\&
Z u= \lim _{t\to +\infty }  e^{ itH_0 } e^{ -it H_\omega
 }  \text{ for any $u\in L^2_c(H_\omega )$} .\tag 2\endalign $$
$W(L^2)=L^2_c(H_\omega )$ is an isomorphism with inverse $Z$.
\endproclaim
{\it Proof.} The existence of $P_c(H_\omega )\circ W$  follows from
Cook's method and Lemma 5.8. By an elementary argument $Wu\in
L^2_c(H_\omega )$ for any $u\in L^2$, so $W=P_c(H_\omega )\circ W$.
We have $W=\widetilde{W} \circ W_1$ with
$$ \align  &
W _1u= \lim _{t\to +\infty }  e^{ itH_q} e^{ -it H_0
 }  u\text{ for any $u\in L^2  (\Bbb R^2 )$} \\& \widetilde{W}  u= \lim _{t\to +\infty }  e^{ itH_\omega \omega} e^{ -it
 H_q
 }  \text{ for any $u\in L^2 $}.\endalign $$
By standard theory $W_1$ is  an isometric isomorphism of $L^2  (\Bbb
R^2 )$ into itself with inverse $Z _1u= \lim _{t\to +\infty }  e^{
itH_0} e^{ -it
 H_q}u$
and by Lemma 5.7 $\widetilde{W}$ is an isomorphism $ L^2  (\Bbb R^2
)\to L^2_c(H_\omega )$ with inverse $\widetilde{Z}$. Then  by
product rule the limit in (2) exists
 and we have
 $Z=Z_1\circ \widetilde{Z}$ with $Z$ the inverse of $W$.

\proclaim{Lemma 5.10} For $u\in L^{2,s}(\Bbb R^2)$ with $s>1/2$ we
have
$$ Wu=u- \frac{1}{2\pi i}\int
_{|\lambda |\ge \omega }  R_{H_\omega } ^{-}(\lambda )  V_\omega
\left [  R_{H_0 } ^{+}(\lambda ) -R_{H_0 } ^{-}(\lambda ) \right ] u
d\lambda.  $$
\endproclaim
{\it Proof.} $Wu \in L^{2 }(\Bbb R^2)$  by Lemma 5.9, but the above
formula is meaningful in the larger space $  L^{2,-s }(\Bbb R^2)$.
For $v\in L^{2,s}(\Bbb
 R^2)\cap L^{2 }_c(H^\ast _\omega ) $ and for $\langle u , v \rangle _{2} =\int _{\Bbb R^2}
  u\cdot \overline{v} dx$ the standard $L^2$
 pairing,  we have  by Plancherel

$$\aligned & \langle W  u,v\rangle _{2}=
\langle  u,v\rangle _{2}+\lim _{\epsilon \to 0^+} \int _0^{+\infty}
\langle V_\omega e^{-iH_0t-\epsilon t } u,
e^{-iH^\ast _\omega t-\epsilon t }v\rangle _{2 }dt\\
=& \langle  u,v\rangle +\lim _{\epsilon \to 0^+ } \frac 1{2\pi }
\int _{-\infty}^{+\infty} \langle AR_{H_0}(  \lambda +i \epsilon )
u, {B} R_{H^\ast _\omega }   (\lambda +i\epsilon ) v\rangle
_{2}d\lambda .
\endaligned$$
By the orthogonality in $L^2(\Bbb R)$ of boundary values of Hardy
functions in $H^2(\Bbb C_+)$ and in $H^2(\Bbb C_-)$ we have for
$\epsilon
>0$
$$\aligned &
\int _{-\infty}^{+\infty}  \langle AR_{H_0}(  \lambda +i \epsilon )
u, {B} R_{H^\ast _\omega }   (\lambda +i\epsilon ) v\rangle _{2}
d\lambda =\\& \int _{-\infty}^{+\infty} \langle  A \left [ R_{H_0}(
\lambda +i \epsilon ) -R_{H_0}(  \lambda -i \epsilon ) \right ]
u,{B} R_{H^\ast _\omega }   (\lambda +i\epsilon ) v\rangle _{2}
d\lambda .
\endaligned
$$
By $u\in L^{2,s}(\Bbb R^2)$  and $v\in L^{2,s}(\Bbb
 R^2)\cap L^{2 }_c(H^\ast _\omega ) $ the limit in the right hand
 side for $ \epsilon \searrow 0$ exists and we have

$$\aligned & \langle W  u,v\rangle _{2}= \langle  u,v\rangle _{2}+\\&  \frac 1{2\pi } \int _{-\infty}^{+\infty}
\langle  A \left [ R_{H_0}( \lambda +i 0) -R_{H_0}(  \lambda -i 0 )
\right ] u,{B} R_{H^\ast _\omega }   (\lambda +i0 ) v\rangle _{2}
d\lambda  =\\& \langle  u,v\rangle _{2}+   \frac 1{2\pi } \int
_{|\lambda |\ge \omega } \langle  A \left [ R_{H_0}( \lambda +i 0)
-R_{H_0}( \lambda -i 0 ) \right ] u,{B} R_{H^\ast _\omega } (\lambda
+i0 ) v\rangle _{2} d\lambda .
\endaligned$$
This yields Lemma 5.10. The crucial part of our linear theory is the
proof  of the following analogue of \cite{Y}:

\proclaim{Lemma 5.11} For any $p\in (1,\infty )$ the restrictions of
$W$ and $Z$ to $L^2\cap L^p$ extend into operators such that for
$C(\omega )<\infty$ semicontinuous in $\omega$ $$\| W\| _{L^p(\Bbb
R^2),L^p_c(H_\omega )}+\| Z\| _{L^p_c(H_\omega ),L^p(\Bbb R^2)
}<C(\omega ).$$ \endproclaim In the next two sections we will
consider $W$ only, since the proof for $Z$ is similar. The argument
in the following two sections is a transposition of \cite{Y}. We
consider diagonal matrices $$E_+=\text{diag}(1 , 0)\text{ and } E_-
=\text{diag}(0 , 1).$$

Keeping in mind Lemma 5.10, $\sigma _1 R (z)=-R(-z )\sigma _1$  for
$R(z)$ equal to  $R_{H_\omega }(z)$ or to $R_{H_0 }(z)$ and
$\sigma_1L^2_c(H_\omega )=L^2_c(H_\omega )$, it is easy to conclude
that the $L^p$ boundness of $W$ is equivalent to $L^p$ boundness of
$$\aligned  Uu:=&\int
_{ \lambda  \ge \omega }  R_{H_\omega } ^{-}(\lambda )  V_\omega
\left [  R_{H_0 } ^{+}(\lambda ) -R_{H_0 } ^{-}(\lambda ) \right ] u
d\lambda  \\= &\int _{ \lambda  \ge \omega }  R_{H_\omega }
^{-}(\lambda )  V_\omega \left [  R_{ 0 } ^{+}(\lambda ) -R_{ 0 }
^{-}(\lambda ) \right ] E_+ u d\lambda .\endaligned$$ As in \cite{Y}
we deal separately with high,treated in \S 6, and low energies,
treated in \S 7. We introduce cut-off functions $\psi_1(x)\in
C_0^\infty(\Bbb R),$ and $\psi_2(x)\in C^\infty(\Bbb R),$ with
$\psi_1(x)+\psi_2(x)=1,$
  $\psi_1(-x)=\psi_1(x),$
$\psi_1(x)=1$ for $|x|\leq C$ and $\psi_1(x)=0$  or $|x|>2C$ for
some $C>\omega $.

\head \S 6   $L^p$ boundness of $U$: high energies
\endhead

This part is almost the same of the corresponding part in \cite{Y2}.
  For $\psi_1(x)$  the cutoff function introduced after Lemma 5.11,   $\psi_1(H_0)$ is a convolution operator with
symbol $\psi_1(|\xi|^2+\omega)$. Both  $\psi_1(H_0)$ and
$\psi_2(H_0)$ are bounded operators in $L^p(\Bbb R^2)$ for any $p\in
[1, \infty].$ In order to estimate the high frequency part (the so
called high energy) $U \psi_2(H_{0}),$ we expand $ R_{H_\omega }
^{-}(\lambda )$ into the sum of few terms of Born series
$$R_{H_\omega } ^{-}(\lambda )=R_{H_0 } ^{-}(\lambda)-R_{H_0 }
^{-}(\lambda)V_\omega R_{H_0 } ^{-}(\lambda)+R_{H_0 }
^{-}(\lambda)V_\omega R_{H_0 } ^{-} (\lambda)V_\omega R_{H_\omega }
^{-}(\lambda ),$$ getting by Lemma 5.10 the decomposition
$U=U_1+U_2+U_3 $ with

$$\aligned
U_1u=- \frac{1}{2\pi i}\int _{\lambda\ge \omega }  R_{H_0}
^{-}(\lambda )  V_\omega    R_{ 0 } ^{+}(\lambda -\omega )   E_+ u
d\lambda,
\endaligned  $$

$$\aligned
U_2u= \frac{1}{2\pi i}\int _{\lambda\ge \omega }  R_{H_0 }
^{-}(\lambda)V_\omega R_{H_0 } ^{-}(\lambda)  V_\omega R_{ 0 }
^{+}(\lambda -\omega ) E_+ u d\lambda,
\endaligned  $$

$$\aligned
& U_3u= - \frac{1}{2\pi i}\int _{\lambda\ge \omega }  R_{H_0 }
^{-}(\lambda)V_\omega R_{H_0 } ^{-}(\lambda)V_\omega R_{H_\omega }
^{-}(\lambda )  V_\omega R_{ 0 } ^{+}(\lambda -\omega )  E_+ u
d\lambda.
\endaligned  $$

 \proclaim{ Lemma 6.1}
The operator $U_1\psi_2(H_{0})$ is bounded in $L^p(\Bbb R^2)$ for
all $1<p<\infty .$ Specifically for any $s>1 $ there exists a
constant $C_s>0$ so that for $T= U_1\psi_2(H_{0})$
$$
\aligned \left\| Tu \right \|_{L^p} \leq C_s\left\| \langle x
\rangle^s V_\omega\right \|_{L^2} \left\| u \right \|_{L^p} \text{
for all $u\in L^p(\Bbb R^2)$}.
\endaligned \tag 1
$$

\endproclaim

{\it Proof.} Recall $R_{0}(z)=(-\Delta -z)^{-1} $ and $
R_{H_0}^{\pm}(z) =      \text{diag}( R_{0}^{\pm}(z-\omega),-
R_{0}^{\pm}(z+\omega)). $ For $u=(u_1,u_2)$,  and for
  $\Cal F$ the Fourier transform, we are reduced to operators of
  schematic form   $\Cal F (E_\pm  {U_1}  u ) (\xi
)=$
$$\aligned
 = \int _ {\lambda \ge \omega } d\lambda \int _{\Bbb R^{2}} \frac 1{
  |\xi |^2+\omega  \mp \lambda  +i0}    \widehat{u}_1 (\xi - \eta )
\delta(\lambda - (|\xi - \eta  |^2+\omega ))
    \widehat{{V}} (\eta ) d\eta ,
\endaligned
$$
with  $  \widehat{{V}} $ the Fourier transform of the generic
component of $V_\omega.$  Then
$$\aligned &
  E_\pm  {U_1}  u    =\int _{\Bbb R^2}d\eta \,   \widehat{{V}} (\eta )
   \, T^{\pm }_{\eta } u_{1\eta }\endaligned $$
where $u_{1\eta }(x)=e^{ix\cdot \eta}u_{1  }(x)$, $T^{- }_{\eta }
u_{1\eta }=\frac{1}{4\pi} K_0(\sqrt{\frac{\eta ^2}{4}+\omega }|\cdot
|)\ast u_{1\eta }$ and by \cite{Y1}
$$T^{+ }_{\eta } u_{1\eta }(x)=\frac{i}{2|\eta |}\int _0^\infty e^{ it|\eta
|} u_{1\eta }(x+t\eta / |\eta |) dt.$$ By \cite{Y2} we have that $T=
E_+  {U_1}$ satisfies inequality (1) while for $T= E_-  {U_1}$ we
use $$\| T^{\pm }_{\eta } u\| _{L^p} \le \frac{1}{4\pi} \left \|
K_0(\sqrt{\frac{\eta ^2}{4}+\omega
 }| x
|) \right \| _{L^1_x} \|  u  _1\| _{L^p}\le C \langle \eta \rangle
^{-1}\|  u  _1\| _{L^p}$$ and so    $\| E_-  {U_1}u \| _{L^p}
\lesssim  \| \widehat{{V}} (\eta )/\langle \eta \rangle \| _{L^1  }
\|  u _1\| _{L^p}.$

\proclaim{ Lemma 6.2} The operator $U_2\psi_2(H_{0})$ is bounded in
$L^p(\Bbb R^2)$ for all $1<p<\infty,$ moreover, there exists a
constant $C_s>0$ so that for $T=U_2\psi_2(H_{0})$
$$
\aligned \left\| T u \right \|_{L^p} \leq C_s\left\| \langle
x \rangle^s V_\omega\right \|^2_{L^2}   \left\| u \right \|_{L^p}\text{ for
all $u\in L^p(\Bbb R^2)$}.
\endaligned \tag 1
$$
is valid, provided $s>1.$
\endproclaim
{\it Proof.}  By \cite{Y1} and with the notation of Lemma 6.1
 we are reduced to a combination of operators

 $$\aligned & I_{\pm , \pm }  u  =\int _{\Bbb R^2}d\eta _1T^{\pm }_{\eta _1}\int _{\Bbb R^2}d\eta _2   \widehat{{V}} (\eta _1)
 \widehat{{V}} (\eta _2-\eta _1)
   \,   T^{\pm }_{\eta _2}u_{1\eta _2 }.\endaligned $$
$Tf=I_{- , - }  u$ satisfies inequality (1) by Proposition 2.2
\cite{Y2} . The other cases follow from Lemma 6.1. For example, for
$K(\eta _1,\eta _2)=\widehat{{V}} (\eta _1)
 \widehat{{V}} (\eta _2-\eta _1)$ and $\widetilde{K}(x,\eta _2)=\int
 d\eta   e^{i\eta \cdot x}K(\eta  ,\eta _2)$,
 $$\aligned & \| I_{\pm , \pm }  u  \| _{L^p}=\| \int _{\Bbb R^2}d\eta _2  \int _{\Bbb R^2}d\eta _1
K(\eta _1 ,\eta _2)  T^{-}_{\eta _1}T^{+}_{\eta _2}u_{1\eta _2 }\|
_{L^p} \\& \le \widehat{C} _s\int _{\Bbb R^2}d\eta _2 \| \langle x
\rangle ^{s}\widetilde{K}(x,\eta _2)\| _{L^2_x} \| T^{+}_{\eta
_2}u_{1\eta _2 }  \| _{L^p} \\& \le \widetilde{C} _s\int _{\Bbb
R^2}d\eta _2 \| \langle x \rangle ^{s}\widetilde{K}(x,\eta _2)\|
_{L^2_x} \langle \eta _2 \rangle ^{-1}   \|  u_{1  }  \|
_{L^p}C_s\left\| \langle x \rangle^s V_\omega\right \|^2_{L^2}
\left\| u _1 \right \|_{L^p}.
\endaligned
$$

\proclaim{ Lemma 6.3}Set $T=U_3\psi_2(H_{0})$. Then $T$ is bounded
in $L^p(\Bbb R^2)$ for all $1\leq p\leq\infty.$
\endproclaim

{\it Proof.}  Schematically
$$\aligned &
  E_+U_3\psi_2(H_{0})u =  \int
_{ k   \ge 0 }  R_{ 0 } ^{-}(k^2) {V}  F (k^2+ \omega )
 {V} \left [  R_{0 } ^{+}(k^2 ) -R_{0 } ^{-}(k^2 ) \right
]\psi_2(\lambda+\omega) u_1 kdk,
\endaligned$$
with $F (k^2+ \omega )=R_{H_0  }^-(k ) {V}  R ^{-}(k )   $ and $V$
the generic component of $V_\omega$. By (3) Lemma 5.4 for $G^{\pm
}_{k,y}(x)= e^{\mp i k|y|}G^{\pm } (x-y,k)$ with $G^{\pm } (x
,k)=\pm \frac{i}{4}H^{\pm}_0(k|x|)$ we have the following analogue
of inequality (3.5) \cite{Y2}
$$ \left | \partial _k^j \langle F(k^2+ \omega )V G^{\pm }_{k,y},VG^{+
}_{k,x}\rangle \right |\le \frac{C_j \| \langle x \rangle ^s
V_\omega \| _{\infty}^3}{k^3 \sqrt{\langle x \rangle \langle y
\rangle}}\tag 1$$ and by Proposition 3.1 \cite{Y2} this yields the
desired result for $T=E_+U_3\psi_2(H_{0}) .$ Since (1) continues to
hold if we replace $G^{+ }_{k,x}$ with  $e^{-ik|x|}G_{k,x}$ with
$G_{k,x}(y)=G(x-y,k), $  where $G(x,k)=
K_{0}(\sqrt{k^2+\omega}|x|)$, we get also the desired result for
$T=E_-U_3\psi_2(H_{0}) .$

\head \S 7   $L^p$ boundness of $U$: Low energies
\endhead

Set $$ Tu:=\int _{ \lambda \ge \omega }  R_{H_\omega } ^{-}(\lambda
) V_\omega \left [  R_{ 0 } ^{+}(\lambda -\omega ) -R_{ 0 }
^{-}(\lambda -\omega ) \right ] \psi_1(\lambda ) E_+u d\lambda   .$$
We want to prove:

 \proclaim{Lemma 7.1} For any $p\in
(1,\infty )$ the restriction of $T$ on $L^2\cap L^p$ extends into an
operator such that $ \| T\| _{L^p(\Bbb R^2),L^p(\Bbb R^2)} <C(\omega
) $ for $C(\omega )<\infty$ semicontinuous in $\omega$.
\endproclaim
 Let $V_\omega =V= \{
V_{\ell j} : \ell,j=1,2 \} $,  $W= \{ W_{\ell j} : \ell,j=1,2 \} $
with $W _{12}=W_{21}=0$,   $W _{22}=1\in \Bbb R$ and $W _{11}(x)=1$
for  $V _{11}(x)\ge 0$  and  $W _{11}(x)=-1$ for  $V _{11}(x)< 0$.
Set $B^\ast =\langle x \rangle   ^{-N}$ for some large $N>0$, and
$A= \{ A_{\ell j} : \ell,j=1,2 \} $ with $A_{11}(x)=|V _{11}(x)|$,
$A_{12}(x)=W _{11}(x) V _{12}(x) $ and $A _{2j}(x)=V _{2j}(x).$ Then
$W^2=1$,  $ B^\ast W A=V$. Let $k>0$ be such that $k^2=\lambda
-\omega $ and set  $M(k)= W+AR_{H_0}^{-}(\lambda )B^\ast $. Then
$$R_{H_\omega}^{-}(\lambda )=R_{H_0}^{-}(\lambda )-R_{H_0}^{-}(\lambda )
B^\ast M^{-1}(k)AR_{H_0}^{-}(\lambda ).$$ We have $ M(k)=W+c^-(k)P+A\widetilde{G_0}B^\ast
+O(k^2\log k)$  where: $c^-(k)=a^-+b^-\log k$;
 $P$ is a projection in $L^2$  defined by
 $$P= \left [ \matrix A_{11}\\
A_{21}\endmatrix \right ] \frac{\langle \cdot , B^{\ast}_{11}\rangle
}{\|
 V_{11} \| _{L^1}};
$$
$$\widetilde{G_0}=\text{diag } \left (  -\frac{1}{2\pi}\log |x|\ast
, - R_0(-2\omega ) \right );$$
$$\| d^j/dk^jO(k^2\log
k)\| _{L^{2 } ,L^{2 } }\le C k ^{2-j}\langle \log k \rangle \quad
j=0,1,2, \quad 0<k<c.$$ Let $Q=1-P$ and let $M_0=W
+A\widetilde{G_0}B^\ast $. Then $QM_0Q$ is invertible in $QL^2$ if
and only if $\omega$ is not a resonance or an eigenvalue for
$H_\omega$ and in that case $ M^{-1}(k)=$
$$  g^{-1}(k) (P-PM_0QD_0Q-QD_0QM_0P  M_0 QD_0Q+QD_0Q+
O(k^2\log k))$$ with $g(k)=c^-\log k+d^-$ for $c^-\neq 0$ and
$D_0=(QM_0Q)^{-1}$ by \cite{JN}. 
We claim now that $QD_0Q-QWQ$ is a Hilbert-Schmidt operator. In fact, following the the argument in Lemma 3 \cite{JY}, we get that the operator $L=P+QM_0Q$ is invertible in $QL^2,$  and $D_0=QL^{-1}Q.$ We have 
$$\aligned & L=W+[A\widetilde{G_0}B^\ast+P+PM_0P-PM_0Q-QM_0P].
\endaligned $$
Set $L:=W(1+\widetilde{S}),$ the operators $P,$ $PM_0P,$ $PM_0Q,$ $QM_0P$ are of rank one while
$A\widetilde{G_0}B^\ast$ is a Hilbert-Schmidt operator. From the fact that $W$ is invertible, we get that also $(1+\widetilde{S})$ is invertible.
Moreover the identity $(1+\widetilde{S})^{-1}=1-\widetilde{S}(1+\widetilde{S})^{-1}$ yields 
$$\aligned & L^{-1}-W=-\widetilde{S}(1+\widetilde{S})^{-1}W
\endaligned ,$$
that is the product of an Hilbert-Schmidt operator with one in $B(L^2(\R^2), L^2(\R^2)).$
Finally, an application of the Theorem VI.22, Chapter VI, in \cite{RS}, shows that $L^{-1}-W$ is of
Hilbert-Schmidt Type.

So we are reduced to
the following list of operators:
$$\aligned & T_0^+u:=\int ^{\infty}_{ 0}  R_{ 0 }
^{-}(k^2  )E_+ V_\omega E_+\left [  R_{ 0 } ^{+}(k^2) -R_{ 0 } ^{-}(k^2 ) \right ] \psi_1(\lambda )
u k dk,
\endaligned $$
and $T_0^-$ defined as above
but with $R_{ 0 }^{-}(k^2  )E_+ $ replaced by $R_{ 0 }(-k^2 -2\omega )E_-$ which are bounded in $L^p$ for $1<p<\infty $ by Lemma 6.1;

$$\aligned &T_1^+u:=\int^{\infty}_{ 0} R_{ 0 }
^{-}(k^2  )E_+ N(k) \left [  R_{ 0 } ^{+}(k^2 ) -R_{ 0 } ^{-}(k^2 )
\right ] \psi_1(\lambda ) E_+u k\, dk
\endaligned $$
with $$\| d^j/dk^jN(k^2\log k)\| _{L^{2, -s } ,L^{2,s } }\le C k
^{2-j}\langle \log k \rangle \quad j=0,1,2, \quad 0<k<c $$ which is
bounded in $L^p$ for $1\le p\le \infty $ by Proposition 4.1
\cite{Y};
$$\aligned &T_2^+u:=\int ^{\infty}_{ 0} R_{ 0 }
^{-}(k^2  )E_+ B^\ast (d(k)F+L+W)A \left [  R_{ 0 } ^{+}(k^2  ) -R_{ 0
} ^{-}(k^2  ) \right ] \psi_1(\lambda ) E_+u k\, dk
\endaligned $$
with $F$ a rank 3 operator, $ L$ a  Hilbert Schmidt operator in
$L^2$, and $d(k)=g^{-1}(k)$. There are also operators $T_j^-,$ for $j=0,1,2,$ defined as above
but with $R_{ 0 }^{-}(k^2  )E_+ $ replaced by $R_{ 0 }(-k^2 -2\omega )E_-$ and bounded in
$L^p.$ So $T_2^{\pm}=T_{2,1}^{\pm}d(\sqrt{-\Delta })+
T_{2,2}^{\pm}+T_{2,3}^{\pm}$  with $T_{2,j}^{\pm}$ for $j=1,2,3$
operators bounded in $L^p$ for $1<p<\infty $ because of the
following statement proved in \cite{Y2} (the $+$ case is exactly
that in \cite{Y2}, and the $-$ case can be proved following the same
argument):

if $K$ is an operator with integral kernel $K(x,y)$ such that for
some $s>1$
$$\| K \| _{s}:=\int _{\Bbb R^2}dy\left (  \int _{\Bbb R^2}dx \langle x \rangle
^{2s}|K(x,x-y)|^2\right ) ^{\frac{1}{2}}<\infty $$ then the operators
$$
\aligned &
Z^+u:=\int _{ 0} ^{\infty}R_{ 0 }
^{-}(k^2  ) K \left [  R_{ 0 } ^{+}(k^2  ) -R_{ 0 } ^{-}(k^2  )
\right ]  u k\, dk \\
& Z^-u:=\int _{ 0} ^{\infty}R_{ 0 }(-k^2+2\omega  ) K \left [  R_{ 0 } ^{+}(k^2  ) -R_{ 0 } ^{-}(k^2  )
\right ]  u k\, dk
\endaligned$$ are bounded in $L^p$ for $1<p<\infty $ with
$\| Z^{\pm} \| _{L^p,L^p}<C_{s,p} \| K \| _{s}.$

\head \S 8   Proofs of Lemmas 3.2, 3.3 and 3.4
\endhead
We mimic Mizumachi \cite{M2}. By the limiting absorption principle
we have

$$
P_c( \omega )e^{-itH_\omega }f =
\frac{1}{2\pi i}\int_{-\infty}^\infty e^{-it\lambda}
(\lambda)P_c( \omega )[R^+_{H_\omega }(\lambda
)-R^-_{H_\omega }(\lambda)]f d\lambda.
  $$
We  consider  a  smooth function  $\chi(x)$  satisfying $0\le
\chi(x)\le 1$  for $x\in\R$, $\chi(x)=1$ if $x\ge 2$ and $\chi(x)=0$
if $x\le 1$.   $\chi_M(x)$ is an even function satisfying
$\chi_M(x)=\chi(x-M)$ for $x\ge0$.  Let
$\widetilde{\chi}_M(x)=1-\chi_M(x).$ We have:

\proclaim{ Lemma 8.1}  For any fixed $s>1$ there exists a positive
  $C(\omega )$ upper semicontinuous in $\omega
, $ such that for any $u\in  {S}(\R^2)$ we have
$$\align &
\| R^\pm_{H_\omega}(\lambda)f\|_{L^2_\lambda(\sigma _c(H_\omega)
;L^{2,-s}_x)} \le C\|f\|_{L^2}
    .
   \endalign $$

\endproclaim
First, we prove Lemma 3.2 assuming Lemma 8.1.

{\it Proof of Lemma 3.2.} We split
$$
P_c( \omega )e^{-itH_\omega }f=P_c( \omega )e^{-itH_\omega }\chi_M(H_\omega
)f+P_c( \omega )e^{-itH_\omega }\widetilde{\chi}_M(H_\omega )f
$$
with
$$\align & P_c( \omega ) \chi_M(H_\omega )e^{-itH_\omega }f=
   \frac{1}{2\pi
i}\int_{-\infty}^\infty e^{-it\lambda} \chi_M(\lambda)(R^+_{H_\omega }(\lambda)-R^-_{H_\omega
}(\lambda))P_c( \omega )f d\lambda,
\\ &
P_c( \omega )e^{-itH_\omega }\widetilde{\chi}_M(H_\omega )f  =
\frac{1}{2\pi i}\int_{-\infty}^\infty e^{-it\lambda}
\widetilde{\chi}_M(\lambda)(R^+_{H_\omega }(\lambda)-R^-_{H_\omega
}(\lambda))P_c( \omega )f d\lambda.
\endalign $$

Integrating  by parts, in $ {S}'_x(\R^2)$ for any $t\ne0$ and $f\in
{S}_x(\R^2 )$
$$
P_c( \omega )e^{-itH_\omega }f=\frac{(it)^{-j}}{2\pi i}
\int_{-\infty}^\infty d\lambda e^{-it\lambda} \pd_\lambda^jP_c(
\omega )\{(R^+_{H_\omega }(\lambda)-R^-_{H_\omega
}(\lambda))\chi_M(\lambda)\}f .
$$  Since by (3) Lemma 5.4 for high energies we have
$$\|\pd_\lambda^jP_c(\omega ) R^\pm_{H_\omega }(\lambda)
: \langle x \rangle ^{(j+1)/2+0} L^{2 }\to \langle x \rangle ^{
-(j+1)/2-0} L^2  ) \| \lesssim \la \lambda\ra^{-(j+1)/2},$$ the above
integral absolutely converges in $\langle x \rangle ^{
-(j+1)/2-0}L^2_x$ for $j\ge2$. Let $g(t,x) \in S(\Bbb R \times\Bbb R^2).$ By
Fubini  and integration by parts,   $j\ge2$,
$$\aligned &
 \la \chi_M(H_\omega ) e^{-itH_\omega }P_c(\omega) f,g\ra _{t,x}
\\=&\frac1{2\pi i}\int_{\R }  dt (it)^{-j} \int_\R d\lambda
e^{-it\lambda} \pd_\lambda^j\left\la \chi_M(\lambda)(R^+_{H_\omega
}(\lambda)-R^-_{H_\omega }(\lambda))f, \overline{g}
\right\ra_x
\\=&
\frac{1}{2\pi i}\int_\R d\lambda  \left\la \pd_\lambda^j   \left \{
\chi_M(\lambda)(R^+_{H_\omega }(\lambda)-R^-_{H_\omega }(\lambda))
  \right \}P_c(\omega) f,\int_\R dt (-it)^{-j}\overline{g }(t) e^{ it\lambda}
\right\ra_x
\\=& \frac{1}{\sqrt{2\pi}i}\int_\R d\lambda
  \left\la
\chi_M(\lambda)(R^+_{H_\omega }(\lambda)-R^-_{H_\omega }(\lambda))P_c(\omega)f,  \overline{\widehat{g}}(\lambda )\right\ra_x.
\endaligned $$
Hence, by Fubini and  Plancherel,  we have

$$\aligned  &
\big |\la\chi_M(H_\omega )
e^{-itH_\omega } P_c(\omega)f,g\ra_{t,x} \big |
 \le
\\\le&(2\pi)^{-1/2}
\|\chi_M(\lambda)(R^+_{H_\omega }(\lambda)-R^-_{H_\omega
}(\lambda))f\|_{L^2_\lambda(\sigma
_c(H_\omega) ;L^{2,-s}_x)} \|
\widehat{g}(\lambda,\cdot) \|_{L^2_\lambda L^{2,s}_x}
\\ = & (2\pi)^{-1/2}
\|\chi_M(\lambda)(R^+_{H_\omega }(\lambda)-R^-_{H_\omega
}(\lambda))f \|_{L^2_\lambda(\sigma
_c(H_\omega) ;L^{2,-s}_x)} \|g\|_{L_t^2L^{2,s}_x},
  \endaligned $$
In a similar way we have
$$\aligned
  &
|\la  e^{-itH_\omega
}\widetilde{\chi}_M(H_\omega)f ,g\ra_{t,x}| \le \\ \le &
(2\pi)^{-1/2}(\|\widetilde{\chi}_M(H_\omega)(R^+_{H_\omega }(\lambda)-R^-_{H_\omega
}(\lambda))f \|_{L^2_\lambda(\sigma
_c(H_\omega) ;L^{2,-s}_x)} \|g\|_{L_t^2L^{2,s}_x},
\endaligned $$
therefore we achieve
$$\aligned
  &  |\la  e^{-itH_\omega }P_c( \omega )f,
 g\ra_{t,x}|\le\\ &\le (2\pi)^{-1/2}\bigl(
\|\chi_M(\lambda)(R_{H_\omega }(\lambda+i0)-R_{H_\omega
}(\lambda-i0))f \|_{L^2_\lambda(\sigma
_c(H_\omega) ;L^{2,-s}_x)}
\\  &
+\|\widetilde{\chi}_M(\lambda)(R^+_{H_\omega }(\lambda)-R^-_{H_\omega
}(\lambda)) f \|_{L^2_\lambda(\sigma
_c(H_\omega) ;L^{2,-s}_x)} \|g\|_{L_t^2L^{2,s}_x}.
\endaligned $$
and  by Lemma 8.1 this estimate yields Lemma 3.2.
\bigskip

{\it Proof of Lemma  3.3} By Plancherel's identity and H\"older
inequalities   we have
$$ \aligned &
\|\int _{0} ^te^{-i(t-s)H_{\omega }}P_c( {\omega
})g(s,\cdot)ds\|_{L^{2,-s}_x L_t^2} \le \\& \le  \| R^+_{H_\omega
}(\lambda)P_c( {\omega })
  \widehat{ \chi }_{[0,+\infty )}\ast _\lambda  \widehat{ g}(\lambda,x)\|_{L^{2,-s}_xL^2_\lambda  }   \le \\& \le  \left\| \,
\|  R^+_{H_\omega }(\lambda)P_c( {\omega }) \|  _{L^{2,s}_x,
L^{2,-s}_x }  \|
   \widehat{ \chi }_{[0,+\infty )}\ast _{\lambda } \widehat{g} (\lambda,x) \|_{L^{2,s}_x}\, \right\|_{L^2_\lambda}.
\endaligned $$
By Lemma 5.4 $ \sup_{\lambda\ge \omega }\|R^+_{H_\omega
}(\lambda)P_c( {\omega })\|_{B(L^{2,s},L^{2,-s})} \lesssim \la
\lambda\ra^{-1/2},
 $  and so
$$ \aligned &
 \sup _{\lambda \in \R }
\| R^+_{H_\omega }(\lambda) P_c( {\omega }) \|_{B(L^{2,s}_x, L^{2,-s}_x)}\| g\|_{L^{2,s}_xL_t^2}
\leq C \| g\|_{L^{2,s}_xL_t^2}.
\endaligned $$

The above inequalities yields Lemma 3.3.

 \medskip

{\it Proof of Lemma  3.4} Let $(q,r)$ be admissible and let $T$ be
an operator defined by
$$Tg(t)=\int_\R ds e^{-i(t-s)H_\omega }P_c( \omega  )g(s).$$
Using Lemmas 3.2 and 3.3 we get
$f:=\int_\R ds e^{isH_\omega }P_c( \omega
)g(s) \in L^2(\R)$  and that there exists a $C>0$
such that
$$
\aligned
\|Tg(t)\|_{L^q_tL^r_x}\le C \|g\|_{L_t^2L_x^{2,s}}
\endaligned \tag 1$$
for every $g\in {S}(\Bbb R\times\Bbb R^2)$.
Since $q>2$, it follows from Lemma 3.1 in [SmS] (see also [Bq]) and (1)
that
$$
\aligned
\left\|\int_{s<t} ds e^{-i(t-s)H_\omega }P_c( \omega
)g(s)\right\|_{L_t^qL_x^p} \lesssim \|g\|_{L_t^2L_x^{2,s}}.
\endaligned
$$
This yields Lemma 3.4 .

\bigskip

To prove Lemma 8.1 observe that it is not restrictive to prove
$$ \| R^\pm_{H_\omega}(\lambda)f\|_{L^2_\lambda((\omega ,\infty
);L^{2,-s}_x)} \le C\|f\|_{L^2}.\tag 8.1$$ Following the argument in
\S 4 \cite{M2} we need the following:

\proclaim{ Lemma 8.2} There exists a positive constant $C$ such that
for $s>1$
$$\aligned
 \|R_{H_0}^\pm(\lambda )f\|_{L^{2,-s}_xL^2_\lambda(\omega ,\infty)} \le
C\|f\|_{L^2}.\endaligned  $$

\endproclaim

{\it Proof.} $E_+R_{H_0}^\pm(\lambda )f= R_{ 0}^\pm(\lambda -\omega
)E_+f$ and by Lemma 4.2 \cite{M2} we get $$ \|
  R_{
0}^\pm(\lambda  )E_+f \|_{L^{2,-s}_xL^2_\lambda(0 ,\infty)}\le C
\sup_x\|R_{ 0}^\pm(\lambda  )E_+f\|_{L^2_\lambda (0 ,\infty)} \le
C\|E_+f\|_{L^2}.\tag 1
$$ We have $E_-R_{H_0}^\pm(\lambda )f= -R_{ 0} (-\omega -\lambda
)E_-f= -\frac{-\Delta +\omega -\lambda }{-\Delta +2\omega +\lambda }
R_{ 0}^+(\lambda -\omega  )  E_ - f  $. So by (1)

$$\aligned & \|
   E_-R_{H_0}^\pm(\lambda )f \|_{L^{2,-s}_xL^2_\lambda(\omega
   ,\infty)}\le \left \|  \frac{-\Delta +\omega -\lambda
}{-\Delta + \omega +\lambda } \right \| _{L^\infty _\lambda ((\omega
,\infty ), B(L^{2,-s}_x ,L^{2,-s}_x))}   \\& \times    \|
  R_{
0}^\pm(\lambda  )E_-f \|_{L^{2,-s}_xL^2_\lambda(0 ,\infty)}\le C_1
\|
  R_{
0}^\pm(\lambda  )E_-f \|_{L^{2,-s}_xL^2_\lambda(0 ,\infty)}\le C_1
C\|E_-f\|_{L^2}.
\endaligned
$$

{\it Proof of inequality (8.1).} We consider the operator
$h_q=-\Delta +q(x)$ introduced in \S 5 and $H_q=\sigma _3(h_q+\omega
)$. We claim that
$$\| R^\pm_{H_q}(\lambda)f\|_{L^2_\lambda ((\omega ,\infty)
,L^{2,-s}_x)} \le C\|f\|_{L^2}.\tag 1$$ Indeed $
E_+R^\pm_{H_q}(\lambda)f=R^\pm_{h_q}(\lambda -\omega )E_+f$  and $
\| R^\pm_{h_q}(\lambda   )E_+f \|_{L^2_\lambda(0,\infty)
,L^{2,-s}_x)} \le C\|f\|_{L^2}$ by Lemma 4.1 \cite{M2}. On the other
hand  $ E_-R^\pm _{H_q}(\lambda)f=$   $$=-R _{h_q}(-\lambda -\omega
)E_-f=-R _{0}(-\lambda -\omega )E_-f+ R _{0}(-\lambda -\omega )qR
_{h_q}(-\lambda -\omega )E_-f.$$ The bound for the first term comes
from Lemma 8.2 and
$$ \aligned & \|  R _{0}(-\lambda -\omega )qR
_{h_q}(-\lambda -\omega )E_-f\| _{L^{2,-s}_x L^2_\lambda   }
\lesssim \|  R _{0}(-\lambda -\omega )qR _{h_q}(-\lambda -\omega
)E_-f\| _{L^{\infty}_x L^2_\lambda   } \\& \lesssim \|   qR
_{h_q}(-\lambda -\omega )E_-f\| _{L^\infty _\lambda  L^{2}_x   } \le
C\| E_-f\| _{   L^{2}_x   }.
\endaligned $$
Armed with inequality (1) we consider   the identity
 $$\aligned & R_{H_\omega } ^{\pm}(\lambda ) =  (1+R_{H_q } ^{\pm}(\lambda
) ( V_\omega -\sigma _3 q)) ^{-1}R_{H_q } ^{\pm}(\lambda )=\\&
=R_{H_q } ^{\pm}(\lambda )-    R_{H_q } ^{\pm}(\lambda )  ( V_\omega
-\sigma _3q )(1+R_{H_q } ^{\pm}(\lambda ) ( V_\omega -\sigma _3 q))
^{-1}R_{H_q } ^{\pm}(\lambda). \endaligned \tag 8.2$$ By (1) it is
enough to bound the     last term in the last sum. This is bounded
by
$$\aligned &  \|  R_{H_q } ^{\pm}(\lambda )  ( V_\omega
-\sigma _3q )(1+R_{H_q } ^{\pm}(\lambda ) ( V_\omega -\sigma _3 q))
^{-1}R_{H_q } ^{\pm}(\lambda) f\| _{L^2_\lambda  L^{2,-s}_x }\le \\&
\|  R_{H_q } ^{\pm}(\lambda )  ( V_\omega -\sigma _3q )(1+R_{H_q }
^{\pm}(\lambda ) ( V_\omega -\sigma _3 q)) ^{-1}\| _{L^\infty
_\lambda B(L^{2,-s}_x,L^{2,-s}_x)} \| R_{H_q } ^{\pm}(\lambda) f\|
_{L^2_\lambda L^{2,-s}_x }   \\& \lesssim \|  R_{H_q }
^{\pm}(\lambda )\| _{L^\infty _\lambda  (  B(L^{2, s}_x,L^{2,-s}_x))
} \| (1+R_{H_q } ^{\pm}(\lambda ) ( V_\omega -\sigma _3 q)) ^{-1}\|
_{L^\infty _\lambda B(L^{2,-s}_x,L^{2,-s}_x)} \|   f\| _{L^2_  x }
\endaligned $$
$\lesssim \|   f\| _{L^2_  x }$ by (1) and by the fact that the
above $L^\infty _\lambda (\omega , \infty )$ norms are bounded by
Lemmas 5.1 and 5.4.

\head  \S 9 Proof of Lemma 4.7 \endhead

The proof is standard and analogous to Lemma 5.8 \cite{Cu2}.
 Recall: \proclaim{Lemma
4.7} We have for $\varphi (x) $ and $ \varphi (t,x)$  Schwarz
functions, for $t\in [0,\infty )$ and for fixed $s>1$ sufficiently
large
$$\aligned &
\|   e^{-iH _\omega t}R_{H_\omega} ^{+}(\Lambda  ) P_c(\omega )
\varphi \| _{L^2_tL^{2,-s}_x} < C(\Lambda , \omega ) \|
  \varphi  (x) \| _{ L^{2, s}_x}   \\&  \left \|
  \int _0^t e^{-iH _\omega (t-\tau )}R_{H_\omega}^{+} (\Lambda  ) P_c(\omega )
\varphi (\tau ) d\tau \right \| _{L^2_tL^{2,-s}_x} < C(\Lambda ,
\omega ) \|
  \varphi (t,x) \| _{L^2_t L^{2, s}_x}  \endaligned
$$
with   $ C(\Lambda , \omega )$ upper semicontinuous in $\omega $ and
in $\Lambda >\omega $.
\endproclaim

{\it Proof.} We consider $\omega < a/<a< <\Lambda <b<\infty$ and
  the partition of unity  $1=g + \widetilde{g} $ with $g\in
C^\infty _0(\Bbb R)$ with $g=1$ in $[a,b]$ and  $g=0$ in $[a/2,2b]$.
By Lemma 3.2 we get
$$\aligned &
\|   e^{-iH _\omega t}R_{H_\omega}^{+} (\Lambda  ) P_c(\omega )
\widetilde{g} (H_\omega )\varphi \| _{L^2_tL^{2,-s}_x}\leq C(\omega)
\| R_{H_\omega}^{+} (\Lambda  ) P_c(\omega ) \widetilde{g} (H_\omega
)\varphi  \| _{L^{2}_x} \\& \le C   (\omega)c_0(a,b,\omega ) \|
\varphi \| _{L^{2}_x}.
\endaligned  $$
Similarly by the proof of Lemma 3.3, for any $s>1$

$$ \aligned &
\|\int _{0} ^te^{-i(t-s)H_{\omega }}R_{H_\omega}^{+} (\Lambda  )P_c(
{\omega })\widetilde{g} (H_\omega )\varphi (s,\cdot)ds\|_{L^{2,-s}_x
L_t^2} \le \\& \le  \| R^+_{H_\omega }(\lambda)R_{H_\omega}^{+}
(\Lambda  )\widetilde{g} (H_\omega )P_c( {\omega })
  \widehat{ \chi }_{[0,+\infty )}\ast _\lambda  \widehat{ \varphi }(\lambda,x)\|_{L^{2,-s}_xL^2_\lambda  }   \le \\& \le  \left\| \,
\|  R^+_{H_\omega }(\lambda)R_{H_\omega}^{+} (\Lambda )\widetilde{g
} (H_\omega )P_c( {\omega }) \|  _{L^{2,s}_x, L^{2,-s}_x }  \|
   \widehat{ \chi }_{[0,+\infty )}\ast _{\lambda } \widehat{\varphi } (\lambda,x) \|_{L^{2,s}_x}\, \right\|_{L^2_\lambda}
   \\& \le
  C(s,a,b,\omega ) \| \varphi \|_{L^{2,s}_xL_t^2}
\endaligned $$
by $(\lambda -\Lambda )R^+_{H_\omega }(\lambda)R_{H_\omega}^{+}
(\Lambda )= R^+_{H_\omega }(\lambda)-R_{H_\omega}^{+} (\Lambda )$,
Lemma 5.4 and
  $|\lambda -\Lambda |\ge a\wedge b$.  We consider now
$$\aligned &\langle x \rangle ^{-\gamma } g(H_\omega )
e^{-iH_\omega t}R_{H_\omega} (  \Lambda +i\epsilon )  P_c(H_\omega )
\langle y \rangle ^{-\gamma }=
\\& e^{-i \Lambda t}\langle x \rangle ^{-\gamma } \int
_t^{+\infty }e^{-i(  H_\omega -\Lambda   -i\epsilon)s}
 g(  H_\omega )
 P_c(H_\omega ) ds\langle y \rangle ^{-\gamma }.
\endaligned \tag 9.1_{\epsilon}
$$
We claim the following:

\proclaim{Lemma 9.1} There are functions  $u(x,\xi )$ defined for
 $x\in \Bbb R^2$ and for $|\xi |\in [a/2,2b]$ with values in $\Bbb C^2$
 such that for any $\chi \in C^\infty _0(a/2,2b)$ we have (for ${^t u
  }\sigma _3f $ the product row column and  $^tu$ the transpose of a column vector)

 $$\chi (H_\omega )f(x)=(2\pi )^{-2} \int _{\Bbb R^4}u (x,\xi )  {^t
 \overline{u }
(y,\xi )}  \sigma _3f(y)
    \chi ( |\xi |^2+\omega )d\xi dy .\tag 9.2$$ There are constants
$c_{\alpha \beta }$ such that
$$|\partial ^\alpha _x
\partial ^\beta _\xi u(x,\xi )| \le c_{\alpha \beta }
\langle x \rangle ^{|\beta |} \quad \text{ for all $x\in \Bbb R^2$
and  $|\xi |\in [a/2,2b]$}.\tag 9.3$$
\endproclaim
Let us assume Lemma 9.1. Then   we can write the kernel of operator
(9.1) as
$$\aligned &\langle x \rangle ^{-\gamma } g(H_\omega )
e^{-iH_\omega t}R_{H_\omega} (  \Lambda +i\epsilon ) \langle y
\rangle ^{-\gamma }=(\text{constant }) \times \\& \langle x \rangle
^{-\gamma } \int _{\Bbb R^3} u(x,\xi )e^{-i(\sigma _3(\xi ^2 +\omega
)-  \Lambda -i\epsilon )s}g(\xi^2+ \omega ) {^t\overline{u}(y,\xi )}
d\xi \langle y \rangle ^{-\gamma }.
\endaligned \tag 9.4
$$
  Estimates (9.3) and elementary integration by parts  yields
$$|(9.4)|\le c\langle x \rangle ^{-\gamma +r}
\langle y \rangle ^{-\gamma +r} s^{-r}e^{-\epsilon t}\, \text{ and
so} \, |(9.1)_{0^+} |\le c\langle x \rangle ^{-\gamma +r} \langle y
\rangle ^{-\gamma +r} \langle t\rangle ^{-r+1}.$$ For $\gamma > r+1$
and $r\ge 3$, we obtain

$$\| e^{-iH_\omega t}R_{H_\omega}^+(  \Lambda  ) g(H_\omega )P_c(H_\omega
)\varphi \| _{L^2_t((0,\infty ), L^{2,-\gamma })}\le C \| \varphi
(x)\| _{  L^{2, \gamma }}.$$ Similarly

$$ \aligned &
\|\int _{0} ^te^{-i(t-s)H_{\omega }}R_{H_\omega}^{+} (\Lambda  )P_c(
{\omega }) {g} (H_\omega )\varphi (s,\cdot)ds\|_{ L_t^2L^{2,-\gamma
}_x} \le \\&  \le \left \|  \int _{0} ^t \langle t-s\rangle ^{-2}\|
\varphi (s,\cdot)ds\|_{ L^{2,\gamma }_x} \right \| _{L^2_t}\le C \|
\varphi \|_{ L_t^2L^{2,\gamma }_x}
\endaligned $$
We need now to prove Lemma 9.1.

\head \S 10 Proof of Lemma 9.1 \endhead

First of all we explain how to define the $u(x,\xi )$. We set
$V_\omega =B^\ast A$ with $A(x)$ and $B^\ast (x)$ rapidly decreasing
and continuous.  Then we have

\proclaim{Lemma 10.1} For any $\lambda > \omega$  and any $ \xi \in
\Bbb R^2$ with   $\lambda =\omega +|\xi | ^2$, in $L^2(\Bbb R^2 )$
the system
$$\left ( 1+A  R _{H_0}^+(\lambda  )B^\ast   \right )
\widetilde{u } = A e^{-i\xi \cdot x}\overrightarrow{e}_1  \tag 1
$$
  admits exactly one solution $\widetilde{u}(x,\xi
)\in  H^2 $ such that for any  $[a,b]\subset (1,\infty ) \setminus
\sigma  _p(H)$ there is a  fixed  $C<\infty $ such that for any
$\lambda \in [a,b]$  and any $\xi$ as above we have
$$ \| \widetilde{u }( \cdot , \xi  ) \|  _{ H^2  }\le C .\tag 2$$

\endproclaim
{\it Proof.} $A  R_{H_0}^+(\lambda  )B^\ast  $ is compact and $\ker
\left ( 1+A  R_{H_0}^+(\lambda  )B^\ast   \right ) =\{ 0\}$ for
$\lambda >\omega$ by  \cite{CPV}, since in that case $\lambda \not
\in \sigma _p(H_\omega )$. By Fredholm alternative we
    get existence and uniqueness of  $\widetilde{u }(x,  \xi  ).$
   Regularity theory and continuity of the coefficients of system (1) with respect to   $\xi$ yield (2)
\bigskip
Let now $^te_1=(1,0)$ and $G_0( |x|
,k)=\text{diag}(\frac{i}{4}H_0^+(k|x|) , -\frac{1}{2\pi}
K_0(\sqrt{k^2+2\omega} |x|))$ for $k>0$. We have $G_0( r ,k)
=\frac{i\sqrt{2}}{4\sqrt{i\pi kr}} e^{ikr}e_1+ O(r^{-\frac{3}{2}})$
and $\partial _rG_0( r ,k)=-k\frac{ \sqrt{2k}}{4\sqrt{i\pi  r}}
e^{ikr}e_1+  O(r^{-\frac{3}{2}})$. We set

$$u  (x,\xi )=e^{-i\xi   \cdot x} {e}_1+v(x,\xi )=e^{-i\xi   \cdot x} {e}_1
-R_{H_0}^+(\lambda  )B^\ast \widetilde{u }( \cdot ,  \xi  ). $$ Then
$(H_\omega -\lambda )u  (x,\xi )=
 B^\ast     \left ( A e^{-i\xi \cdot x} {e}_1-\widetilde{u }
 -A  R_{H_0}^{+}(\lambda  )B^\ast    \widetilde{u }
\right )=0. $ Notice $B^\ast \widetilde{u }  =V_\omega u$ so $v(x,
\xi )=e^{-ix\cdot \xi}w(x, \xi )$ where $w(x, \xi )$ is the unique
solution in $L^2_{-s}$, $s> 1 $, of the integral equation
$$\aligned & w(x,\xi )= -F(x,\xi )-
\int _{\Bbb R^2}
 G_0( |x-z|,|\xi |)
e^{i(x-z)\cdot \xi }V _{\omega} (z) w(z,\xi ) dz,
\endaligned \tag 1
$$
with $$F(x,\xi )= \int _{\Bbb R^2} G_0( |x-z|,|\xi |)V_{\omega }(z)
e^{i(x-z)\cdot \xi }e_1 dz.$$   It is elementary to show that, for
$|\xi |\in [a,b]$, then $|\partial ^\alpha _x
\partial ^\beta _\xi F(x,\xi )| \le \tilde c_{\alpha \beta }
\langle x \rangle ^{|\beta |-1/2}.$ By standard arguments and Lemmas
5.3 and 5.4 we have $|\partial ^\alpha _x
\partial ^\beta _\xi w(x,\xi )| \le \tilde c_{\alpha \beta }
\langle x \rangle ^{|\beta |}.$ This yields (9.3). To get (9.2) we
follow the presentation  in Chapter 9 \cite{Ta}. We denote by
$R_{H_\omega }^{\pm }(x,y,k  )$  the kernel of $ R_{H_\omega }^{\pm
}(k^2+\omega  )$. We set
$$R_{H_\omega }^{+}(x,y,k  )=G_0(|x-y|,k  )+h(x,y,k  )$$ with
$  h(\cdot ,y,k  ) =-R^{+}_{H_0}(k^2+\omega )V_\omega G_0(|\cdot
-y|,k  ).$ Let $(r, \Sigma)$ be polar coordinates on the sphere
$S^1,$ then we claim:

\proclaim {Lemma 10.2} Let $k>0$. For $r\to \infty$ we have uniform
convergence on compact sets of, with $u\cdot (1,0)$ the  raw column
product between column $u$ and raw $(1,0)$,
$$\align & R_{H_\omega }^{+}(x,r\Sigma ,k )= \frac{i\sqrt{2}}{4\sqrt{i\pi kr}}e^{ikr}u(x,k\Sigma  )\cdot (1,0)+O(r^{-2})\tag 1\\&
\frac{\partial}{\partial r }R_{H_\omega }^{+}(x,r\Sigma ,k )=-\frac{
\sqrt{2}}{4\sqrt{i\pi kr}} ke^{ikr}u(x,k\Sigma )\cdot
(1,0)+O(r^{-2}),\tag 2 \\& R_{H_\omega }^{+}( r\Sigma , y ,k )=
\frac{i\sqrt{2}}{4\sqrt{i\pi kr}}e^{ikr} \left [ \matrix 1
\\ 0 \endmatrix \right ] {^tu}(y,k\Sigma ) \sigma _3 +O(r^{-2}), \tag 3
\\&\frac{\partial}{\partial r }R_{H_\omega }^{+}( r\Sigma , y ,k ) =-\frac{
\sqrt{2}}{4\sqrt{i\pi kr}} ke^{ikr}\left [ \matrix 1
\\ 0 \endmatrix \right ] {^tu} (y,k\Sigma ) \sigma _3 +O(r^{-2}).\tag 4
\endalign
$$
For $R_{H_\omega }^{-}(x,y,k  )$ the asymptotic expansion follows
from   $R_{H_\omega }^{-}(x,y,k  )=\overline{R _{H_\omega
}^{+}}(x,y,k  ).$
\endproclaim

We write $R_{H_\omega }^{+}(x,r\Sigma ,k )=G_0(|x-r\Sigma |,k
 )+h(x,r\Sigma ,k )$ with
$$\aligned & h(x,r\Sigma ,k  )  =-R^{+}_{H_0}(k^2+\omega )V_\omega  G_0(|\cdot -r\Sigma |,k
 )\\& =-R^{+}_{H_0}(k^2+\omega )\left
[ V_\omega (x) \left (\frac{i\sqrt{2}}{4\sqrt{i\pi kr}}
e^{ikr}e^{-ik\Sigma \cdot x}\text{diag}(1,0)+ O(r^{-\frac{3}{2}})
\right ) \right ] .\endaligned$$ We have $$\| V_\omega
(x)G_0(|x-r\Sigma |,k
 )-V_\omega (x)  \frac{i\sqrt{2}}{4\sqrt{i\pi kr}} e^{ikr} e^{-ik\Sigma \cdot x}
 \text{diag}(1,0)\| _{L ^{2,s}_{x} }=O(r^{-3/2}).$$
From $ v(x,\xi )= -R^{+}_{H_0}(k^2+\omega ) V_{\omega
}(x)e^{-ik\Sigma \cdot x}e_1 $, with $ {^te_1}=(1,0) $  we get
$v(x,\xi ) \, {^te_1}= -R^{+}_{H_0}(k^2+\omega ) V_{\omega
}(x)e^{-ik\Sigma \cdot x}\text{diag}(1,0).$ Then we
 conclude for any $s>1$

$$\| h(x,r\Sigma ,k   )- \frac{i\sqrt{2}}{4\sqrt{i\pi kr}}v(x,k\Sigma  ) {^te_1}\| _{L ^{2,-s} }=O(r^{-3/2})$$
and
$$\| R_{H_\omega }^{+}(x,r\Sigma ,k   )
-\frac{i\sqrt{2}}{4\sqrt{i\pi kr}}u(x,k\Sigma  ) {^te_1}\| _{L ^{2,-s} }=O(r^{-3/2}).$$ Then
    point wise
$h(x,r\Sigma ,k +i0 )- \frac{i\sqrt{2}}{4\sqrt{i\pi kr}}v(x,k\Sigma
){^te_1} =O(r^{-3/2})$ and
$$R_{H_\omega }^{+}(x,r\Sigma ,k   )-\frac{i\sqrt{2}}{4\sqrt{i\pi kr}}u(x,k\Sigma  ) {^te_1} =O(r^{-3/2}).$$
This yields (1) in Lemma 10.2. (2) can be obtained with a similar
argument. (3) and (4) follow from (1) and (2) by
$$\sigma _3R_{H_\omega }^{\pm }(x,y ,k   )\sigma _3= R_{H_\omega ^\ast}^{\pm }(x,y ,k
)= {^tR} _{H_\omega }^{\mp }(y,x ,k   ).$$

By Lemma 3.5 for $v\in L^2(H_\omega ) \cap C^\infty _0$ and for
$\varphi \in C^\infty _0(\Bbb R)$ supported in $(\omega ,\infty)$ we
have
$$\aligned \varphi (H_\omega )v(x) = \frac{2}{\pi }\int _0^\infty k \,
dk\int _{\Bbb R^2}\varphi (k^2+\omega )\Im R_{H_\omega }^{+ }(x,y,k
) v(y)dy.
\endaligned $$
We prove (here $u \, {^t{\overline{u}}}$ is a  raw column product
between column $u$ and raw $^t{\overline{u}}$)
$$\aligned \Im R_{H_\omega }^{+ }(x,y,k  )=\frac{1}{8\pi }\int _{S^1}u(x,k\Sigma  )\, {^t{\overline{u}(y,k\Sigma
) }} \sigma _3d\Sigma,  \endaligned \tag 3$$ where $d\Sigma $ is the
standard measure on $S^1.$ By the Green theorem   for $S_R=\{ z\in
\Bbb R^2: |z|=R\}$, $|x|<R$, $|y|<R$ and $r=|z|$

By Green theorem   for $S_R=\{ z\in \Bbb
R^2: |z|=R\}$, $|x|<R$ and $|y|<R$,

$$\aligned &\Im R_{H_\omega }^{+ }(x,y,k )= \frac{1}{2i}  \int _{  S_R}
I(x,y,z,k)  d\ell (z)\\& I(x,y,z,k):=  R_{H_\omega }^{+ }(x,z,k )
\sigma _3\partial _{|z|}
  R_{H_\omega }^{- }(z,y,k )   -(
\partial _{|z|} R_{H_\omega }^{+ }(x,z,k ) )\sigma _3
R_{H_\omega }^{- }( z,y,k )
\endaligned $$

By Lemma 10.2

$$\aligned &
\left |\Im R_{H_\omega }^{+ }(x,y,k )-\frac{1}{8\pi}\int _{S^1}u(x,k\Sigma  )\,
{^t{\overline{u}(y,k\Sigma ) }}\sigma _3d\Sigma \right |=\\& = \left
| \frac{R}{2i} \int _{  S^1} I(x,y,r\Sigma,k)|_{r=R} d\Sigma
-\frac{1}{8\pi}\int _{S^1}u(x,k\Sigma  )\,
{^t{\overline{u}(y,k\Sigma ) }}\sigma _3 d\Sigma \right |\leq
O(R^{-\frac{3}{2}}).
\endaligned $$
Therefore, taking $R\rightarrow +\infty,$ we arrive at (3).
Moreover, we obtain
$$\aligned &
 \varphi (H_\omega )v(x) = \frac{2}{\pi }\int _0^\infty k \,
dk\int _{\Bbb R^2}\varphi (k^2+\omega )\Im G(x,y,k  ) v(y)dy=\\& =
\frac{1}{4\pi^2 }\int _0^\infty k \, dk\int _{\Bbb R^2} \int
_{S^1}u(x,k\Sigma  )\, {^t{\overline{u}(y,k\Sigma ) }}  \sigma
_3v(y) \varphi (k^2+\omega) d\Sigma dy=\\& =(2\pi )^{-2} \int _{\Bbb
R^4}u (x,\xi )  {^t
 \overline{u }
(y,\xi )}  \sigma _3v(y)
    \varphi ( |\xi |^2+\omega )d\xi dy,
\endaligned $$
that is the integral representation (9.2).
 This completes the proof
of   Lemma 9.1.

\enddocument